\documentclass{article}

\usepackage{arxiv}

\usepackage[utf8]{inputenc} 
\usepackage[T1]{fontenc}    
\usepackage{hyperref}       
\usepackage{url}            
\usepackage{booktabs}       
\usepackage{amsfonts}       
\usepackage{nicefrac}       
\usepackage{microtype}      
\usepackage{lipsum}		
\usepackage{amsthm}
\usepackage{graphicx}
\usepackage{doi}
\usepackage{physics}
\usepackage{xcolor}
\usepackage{algorithm}
\usepackage{algpseudocode}
\usepackage{amssymb}
\usepackage{pifont}
\newcommand{\pd}[2]{\frac{\partial #1}{\partial #2}}
\newcommand{\pdd}[2]{\frac{\partial^2 #1}{\partial #2^2}}
\newcommand{\cmark}{\ding{51}}%
\newcommand{\xmark}{\ding{55}}%

\usepackage{setspace}
\title{Boundary condition enforcement with PINNs: a survey, critical comparison, and verification on 3D geometries}

\date{} 					

\author{
    Conor Rowan \\
	Smead Aerospace Engineering Sciences\\
    3775 Discovery Drive \\
	University of Colorado Boulder\\
	Boulder, CO 80309 \\
	\texttt{conor.rowan@colorado.edu} \\
     \\
    \And
    Kai Hampleman \\
	Smead Aerospace Engineering Sciences\\
    3775 Discovery Drive \\
	University of Colorado Boulder\\
	Boulder, CO 80309 \\
	\texttt{kai.hampleman@colorado.edu}
    \And
    Kurt Maute \\
	Smead Aerospace Engineering Sciences\\
    3775 Discovery Drive \\
	University of Colorado Boulder\\
	Boulder, CO 80309 \\
	\texttt{kurt.maute@colorado.edu} \\
    \And
    Alireza Doostan \\
	Smead Aerospace Engineering Sciences\\
    3775 Discovery Drive \\
	University of Colorado Boulder\\
	Boulder, CO 80309 \\ \texttt{alireza.doostan@colorado.edu}
     \\
}

\renewcommand{\headeright}{}
\renewcommand{\undertitle}{}
\renewcommand{\shorttitle}{}

\hypersetup{
pdftitle={A template for the arxiv style},
pdfsubject={q-bio.NC, q-bio.QM},
pdfauthor={David S.~Hippocampus, Elias D.~Striatum},
pdfkeywords={First keyword, Second keyword, More},
}

\begin{document}
\maketitle

\begin{abstract}
    Since their advent nearly a decade ago, physics-informed neural networks (PINNs) have been studied extensively as a novel technique for solving forward and inverse problems in physics and engineering. The neural network discretization of the solution field is naturally adaptive and avoids meshing the computational domain, which can both improve the accuracy of the numerical solution and streamline implementation. However, there have been limited studies of PINNs on complex three-dimensional geometries, as the lack of mesh and the reliance on the strong form of the partial differential equation (PDE) make boundary condition (BC) enforcement challenging. Techniques to enforce BCs with PINNs have proliferated in the literature, but a critical discussion of these techniques and a study of their efficacy on geometrically complex three-dimensional test problems are lacking. In this work, we i) comprehensively survey BC enforcement strategies for PINNs and perform a comparison on a relevant subset, ii) propose a general solution framework for arbitrary three-dimensional geometries based on weak boundary enforcement, and iii) verify the methodology on three-dimensional, linear and nonlinear test problems with combinations of Dirichlet, Neumann, and Robin boundaries. Our approach is agnostic to the underlying PDE, the geometry of the computational domain, and the nature of the BCs, while requiring minimal hyperparameter tuning. This work represents a step in the direction of establishing PINNs as a mature numerical method, capable of competing head-to-head with incumbents such as the finite element method. 
\end{abstract}

\keywords{Physics-informed neural networks \and Constrained optimization \and Boundary condition enforcement \and Numerical methods}


\section{Introduction}
\label{sec:intro}

\paragraph{} Neural networks were first introduced more than half a century ago as flexible function approximators for data-driven classification and regression tasks \cite{rosenblatt_perceptron_1958}. In the past decade, physics-informed neural networks (PINNs) have emerged from the scientific machine learning community as a modification to the standard machine learning workflow. PINNs use physical laws in the form of governing partial differential equations (PDEs) in the process of training a neural network to represent the solution field. Consistent with their data-driven origins, PINNs often supplement the physics-based objective with a second term that enforces agreement with available measurement data. An example of this is the field of inverse problems, where trainable parameters are introduced into the governing equation and estimated by simultaneously minimizing the PDE residual and the mismatch of the solution with measurement data. The efficacy of this method for inverse problems was first demonstrated for estimating the unknown viscosity of the Navier-Stokes equation \cite{raissi_physics-informed_2019}. It has since been extended to conjugate heat transfer \cite{cai_physics-informed_2021}, contact mechanics \cite{sahin_solving_2024}, and beam bending \cite{zhou_data-guided_2024}, to name a few. 


\paragraph{} While PINNs are often used to blend physics and data, they can also be used as a numerical method analogous to finite element, finite volume, finite difference, or spectral methods. In this setting, no data is used to formulate the objective/loss function, and the neural network simply acts as a discretization of the solution field, with parameters tuned to satisfy the governing PDE. Lagaris et al. performed prescient early studies of neural network discretizations of partial differential equations (PDEs) over two decades ago \cite{lagaris_artificial_1998, lagaris_neural-network_2000}. In their work, they discretize the PDE solution with a multilayer perceptron (MLP) neural network and determine its weights and biases by minimizing the sum of the squares of the PDE residual at discrete collocation points. More than 20 years later, the term PINNs was coined when, with the help of modern machine learning libraries, these methods were extended to more complex physical models \cite{raissi_physics-informed_2019}. We note that the “Deep Galerkin method” employs similar strategies and was introduced around the same time \cite{sirignano_dgm_2018}. Follow-up works have extended these methods to different problems, such as elasticity \cite{cai_deep_2020, kag_physics-informed_2024}, contact mechanics \cite{sahin_solving_2024}, heat transfer \cite{madir_physics_2024}, fluid mechanics \cite{jin_nsfnets_2021}, homogenization \cite{leung_nh-pinn_2022}, and structural reliability analysis \cite{meng_pinn-form_2023}. Owing to concerns about the scalability of the vanilla PINNs approach, domain decomposition methods were introduced and explored in \cite{ameya_d_jagtap_extended_2020, hu_when_2022, jagtap_conservative_2020}. Various other modifications to the network architecture and loss function have been proposed to improve the accuracy and convergence rates of the learned solution \cite{anagnostopoulos_residual-based_2024, wang_respecting_2024, li_physical_2023, gao_physics-informed_2022, du_evolutional_2021}.

\paragraph{} For many physical systems of engineering interest, the solution of the differential equation corresponds to a minimum of an energy functional. Noting that this energy functional has lower-order spatial derivatives than the corresponding strong form loss---and thus reduces the cost of forward evaluations and backpropagation--the energy defines an attractive objective function for solutions discretized with neural networks. This technique was first introduced as the “Deep Ritz Method” \cite{e_deep_2017} and has been studied in the context of linear elasticity \cite{liu_deep_2023}, hyperelasticity \cite{abueidda_deep_2022, nguyen-thanh_deep_2020}, thermoelasticity \cite{lin_investigating_2024}, and fracture mechanics \cite{manav_phase-field_2024, ghaffari_motlagh_deep_2023}. We note that such energy functionals do not exist for all physical systems, including but not limited to fluid systems and solids with dissipative constitutive models.

\paragraph{} Beyond the strong form and energy loss functions, neural networks can also be trained with the weak form of the governing PDE. Instead of a scalar loss function, the weak form is a system of equations expressing the orthogonality of the PDE residual to a chosen basis of ``test'' functions. As for the energy formulation, the weak form benefits from lower-order spatial differentiation, yet, because it is derived from the strong form of the governing equation, it exists for all physical systems. Inspired by adversarial training from computer vision, an interesting formulation of the weak form is given in \cite{zang_weak_2020}, where the test function is a neural network trained to maximize the weighted residual of the PDE. In variational physics-informed neural networks (VPINNs), the spacetime solution is parameterized as a neural network, the strong form is integrated against a suitable basis of test functions, and the neural network parameters are found by minimizing the norm of the weak form residual \cite{kharazmi_variational_2019}. Follow-up works have explored the Petrov-Galerkin weak form loss in the context of a variety of other problems from engineering mechanics \cite{kharazmi_hp-vpinns_2021, shang_deep_2022, shang_randomized_2023, khodayi-mehr_varnet_2019}.


\paragraph{} Enforcing boundary conditions (BCs) has been one of the primary obstacles to using PINNs for PDEs defined on complex two- and three-dimensional geometries. With the finite element method (FEM), Dirichlet BCs are built into the discretization and Neumann BCs are enforced weakly, meaning that they appear as forcing terms in the governing equation for equilibrium \cite{hughes_thomas_j_r_finite_2000}. Given that a PINN solution is mesh-free, it is not straightforward to prescribe the solution along the Dirichlet region of the boundary. When using the variational energy or weak form as an objective, the Neumann BC appears naturally in the objective function, but this is not the case for the strong form loss. Thus, the PINN methodology leaves open important questions about how to best handle BCs. In the original PINN formulation, both Dirichlet and Neumann boundaries are enforced with penalties added to the strong form loss \cite{raissi_physics-informed_2019}. The penalty approach has been adopted by several other authors for both static and dynamic problems \cite{leung_nh-pinn_2022, meng_pinn-form_2023, jin_nsfnets_2021}. Noting that multi-objective optimization problems can lead to stiff gradient flow dynamics, various strategies have been proposed to adaptively weight the different terms in the loss function \cite{bischof_multi-objective_2025, chen_gradnorm_2018, zhou_dual-balancing_2025}. Of particular relevance, Wang et al. \cite{wang_understanding_2020} proposed a simple algorithm that dynamically updates penalty weights so that the gradient contributions of each term remain balanced.

\paragraph{} Since BCs can be cast as constraints on the solution field, Lagrange multipliers can be used as another enforcement strategy, though this gives the objective function a saddle point structure \cite{ruszczynski_nonlinear_2006}. A ``Deep Uzawa algorithm'' is proposed in \cite{makridakis_deep_2024}, whereby a stationary point of the Lagrange function objective is found through a staggered descent-ascent strategy. A method that blends penalty methods and Lagrange multipliers is the ``Augmented Lagrangian'' method, which enforces constraints by repeatedly solving modified optimization problems in an iterative scheme \cite{hu_conditionally_2025, basir_physics_2022, son_enhanced_2023, lu_physics-informed_2021, zhang_physics-informed_2025}. Though the constraints are enforced accurately, a downside of the Augmented Lagrangian is that the neural network must be repeatedly trained at different settings of the penalty parameter and Lagrange multiplier. A modification to the Augmented Lagrangian for weak formulations of the PDE solution is ``Nitsche's Method,'' which was investigated for high-dimensional and nonlinear problems in \cite{liao_deep_2021}. Though Nitsche's method is intended to improve the convergence properties of standard penalty formulations, the authors in \cite{berrone_enforcing_2023} reported that it did not have this effect. Furthermore, Nitsche's method is limited in scope as it only applies to weak and variational formulations of the PDE system. Another related strategy, which combines aspects of penalty methods and Lagrange multipliers, is ``Self-adaptive PINNs'' (SA-PINNs), which introduces trainable penalty parameters that increase until the BCs are satisfied \cite{mcclenny_self-adaptive_2023}.

\paragraph{} All of the above methods rely on modifying the objective function and/or optimization procedure such that the stationary point which the optimizer targets corresponds to a solution that exactly or approximately satisfies the BCs. An alternative approach is to modify the discretization such that the solution field satisfies the BCs automatically. This is possible in the case of Dirichlet BCs, where the neural network discretization of the solution can be multiplied by a ``distance function'' which is zero along the boundary and positive inside the domain \cite{wang_exact_2023, sukumar_exact_2022, sheng_pfnn_2021, berrone_enforcing_2023}. This enforces homogeneous Dirichlet BCs. For inhomogeneous Dirichlet boundaries, a second function can be added into the discretization after multiplying the neural network by a distance function, thus ensuring that the BCs are satisfied for any setting of the neural network parameters. Though this approach leads to rapid convergence and accurate PDE solutions, it may be non-trivial to construct the distance function and the additive term for the inhomogeneous boundary. Furthermore, it does not extend to Neumann boundaries and becomes especially cumbersome for domains combining Dirichlet and other boundary types. Going forward, we refer to domains with combinations of Dirichlet, Neumann, and/or Robin boundaries as having ``mixed'' BCs. To the best of the authors' knowledge, using distance functions in this setting of mixed boundaries has not been discussed in the PINNs literature.
 

\paragraph{} Determining best practices for enforcing BCs on complex geometries is a crucial step toward making PINNs competitive with traditional numerical methods such as FEM. Accordingly, the focus of this work will be on studying and improving the capabilities of PINNs as a numerical method with an emphasis on BC enforcement. Our review suggests that a number of formulations of the physics loss and BC enforcement have been proposed, but the scientific machine learning community has not come to a consensus on best practices. Specifically, our goal is to move in the direction of a general-purpose PINN-based numerical solver, which can handle a wide range of governing equations, geometries, and BCs with minimal hyperparameter tuning. To do this, we must address the three fundamental tasks of i) formulating the physics loss, ii) defining the geometry, and iii) enforcing BCs. To this end, our contributions in this paper are as follows:

\begin{enumerate}
    \item We critically review the three primary formulations of the physics loss (strong form, weak form, variational energy) and suggest that for a PINN-based solver to be generally applicable and to take full advantage of the neural network discretization, it must be based on the strong form loss;
    \item We comprehensively survey techniques available for BC enforcement and perform a side-by-side performance comparison of a relevant subset of these techniques on a canonical test problem;
    \item Informed by the results of this study, we identify best practices for BC enforcement in the context of PINNs and, by bringing together existing methods in the literature, provide a framework to solve PDEs on arbitrary three-dimensional geometries;
    \item Relying only on a small set of clearly-defined hyperparameters, we verify this methodology on three linear and nonlinear three-dimensional problems from engineering mechanics with complex geometries and BC configurations.
\end{enumerate}

The rest of this paper is organized as follows. In Section \ref{sec:strong form}, we justify our choice to use the strong form loss as the objective function, arguing that the applicability of the Deep Ritz Method is too limited and that the weak form undermines many of the advantages of neural network discretizations. In Section \ref{sec:comparison}, we discuss the various techniques for BC enforcement and perform a side-by-side comparison of a relevant subset of these techniques on a two-dimensional test problem. In Section \ref{sec:method}, we build on the results of the comparison to illustrate a general-purpose solution methodology based on PINNs. The few hyperparameters of the solution methodology will be clearly identified and discussed. In Section \ref{sec:examples}, we verify the method against manufactured solutions on three-dimensional example problems, as well as compare its performance to a second technique for boundary condition enforcement. Finally, in Section \ref{sec:conclusion}, we close with concluding remarks and directions for future work.


\section{On the strong form loss}
\label{sec:strong form}

\paragraph{} As discussed above, several formulations of the physics loss have been explored in the PINNs community. All of these formulations are instances or variants of the strong form, weak form, or energy objective. We argue that the strong form loss provides the most general and flexible foundation for a PINN-based numerical solver intended to handle a broad class of PDEs. We reiterate that the energy objective only exists for some physical systems (steady-state heat transfer, linear and hyperelasticity, variational fracture mechanics, etc.). Notably, even for dynamics problems with variational principles, the stationary point of the energy functional is often a saddle point, not a minimum. The so-called ``principle of least action'' from classical mechanics is indeed a principle of stationary action \cite{goldstein_classical_2002}. Worse, the structure of the saddle is not known, meaning that, unlike saddle problems from Lagrange multipliers, it is not clear over which variables to maximize and over which to minimize. As such, the analyst is prevented from finding a stationary point with a min-max formulation of the optimization problem. This explains why the Deep Ritz method, which relies on \textit{minimizing} an energy functional, has not been used on time-dependent systems. Thus, not only is Deep Ritz infeasible for fluid systems and non-conservative solid constitutive models, but it also cannot be used as a general strategy for dynamical systems. For any PINN-based solver that aspires to maximize its scope, an energy formulation must be ruled out. 

\paragraph{} Though the weak form exists for all PDEs, the quality of the solution relies heavily on the choice of test functions. The Bubnov-Galerkin weak form asserts that the PDE residual is orthogonal to the local tangent of the approximation space. For linear discretizations, this reduces to enforcing orthogonality of the residual against each basis function. In nonlinear discretizations---where parameters define both the basis functions and their coefficients---this condition can admit trivial solutions \cite{rowan_physics-informed_2025}. In such cases, the residual may be reduced not by satisfying the PDE but by learning poor basis functions that weaken the approximation. Knowingly or not, this has pushed researchers to rely on the Petrov-Galerkin weak form, which makes use of distinct approximation and test spaces. One problem with this approach is that it requires the offline construction of a basis to test the PDE residual against. In \cite{shang_deep_2022, khodayi-mehr_varnet_2019}, the test functions are taken to be piecewise polynomials defined on a finite element mesh, which arguably voids the ``mesh-free'' promise of PINNs. Mesh-free spectral bases can be used as test functions, but such orthogonal functions are only tabulated for simple domain geometries \cite{kharazmi_hp-vpinns_2021, berrone_variational_2022}. One work has investigated using ``geometry-informed'' neural network discretizations, but this method relies on offline FEM computations to obtain a spectral basis from eigenfunctions of the Laplace operator, which again calls into question whether the method is truly mesh-free \cite{sahli_costabal_pinns_2024}. 

\paragraph{} An even more fundamental problem with the Petrov-Galerkin weak form is that it may not take full advantage of the expressivity of the neural network discretization. To see this, consider the following elliptic boundary value problem:
\begin{equation}\label{bar}
    \pd{}{x}\qty( \qty(1 + \frac{1}{2} \sin(20\pi x)) \pd{u}{x} ) + 100 \sin( \pi x) = 0, \quad x \in[0,1], \quad u(0)=u(1)=0.
\end{equation}
This is the governing equation for the displacement field $u(x)$ of an axially loaded bar of unit length made of a material with multiscale properties. An exact solution can be obtained by directly integrating Eq. \eqref{bar} and choosing integration constants to enforce the homogeneous Dirichlet BCs. Now, we obtain a PINN-based solution with both strong and weak form objectives. In both cases, the displacement is discretized with
\begin{equation}\label{discret}
    \hat u(x;\boldsymbol \theta) = \sin( \pi x) \mathcal{N}(x;\boldsymbol \theta),
\end{equation}
\noindent where $\mathcal{N}(x;\boldsymbol \theta)$ is an MLP neural network with parameters $\boldsymbol \theta$. The discretization of Eq. \eqref{discret} satisfies the Dirichlet boundaries automatically. The strong form loss is given by 
\begin{equation*}
    \mathcal{L}^{\text{STRONG}}(\boldsymbol \theta) = \int_0^1 \qty[\pd{}{x}\qty( \qty(1 + \frac{1}{2} \sin(20\pi x)) \pd{\hat u}{x} ) + 100 \sin( \pi x)]^2 dx.
\end{equation*}

Defining a set of test functions $\{ v_i(x)\}_{i=1}^N$, the weak form loss is the squared magnitude of the residual of the Petrov-Galerkin orthogonality condition, given by
\begin{equation*}
    \mathcal{L}^{\text{WEAK}}(\boldsymbol \theta) = \frac{1}{2}\sum_{i=1}^N \qty( \int_0^1 \qty[\pd{}{x}\qty( \qty(1 + \frac{1}{2} \sin(20\pi x)) \pd{\hat u}{x} ) + 100 \sin( \pi x)] v_i(x) dx)^2.
\end{equation*}
The test functions $v_i(x)$ are defined as piecewise linear finite element hat functions corresponding to a solution with homogeneous Dirichlet boundaries:
\begin{equation*}
    v_i(x) = \text{max}\Big( 0 , 1-\Big|x-\frac{i}{N+1}\Big|(N+1)\Big).
\end{equation*}

\begin{figure}[hbt!]
\centering
\includegraphics[width=1.0\textwidth]{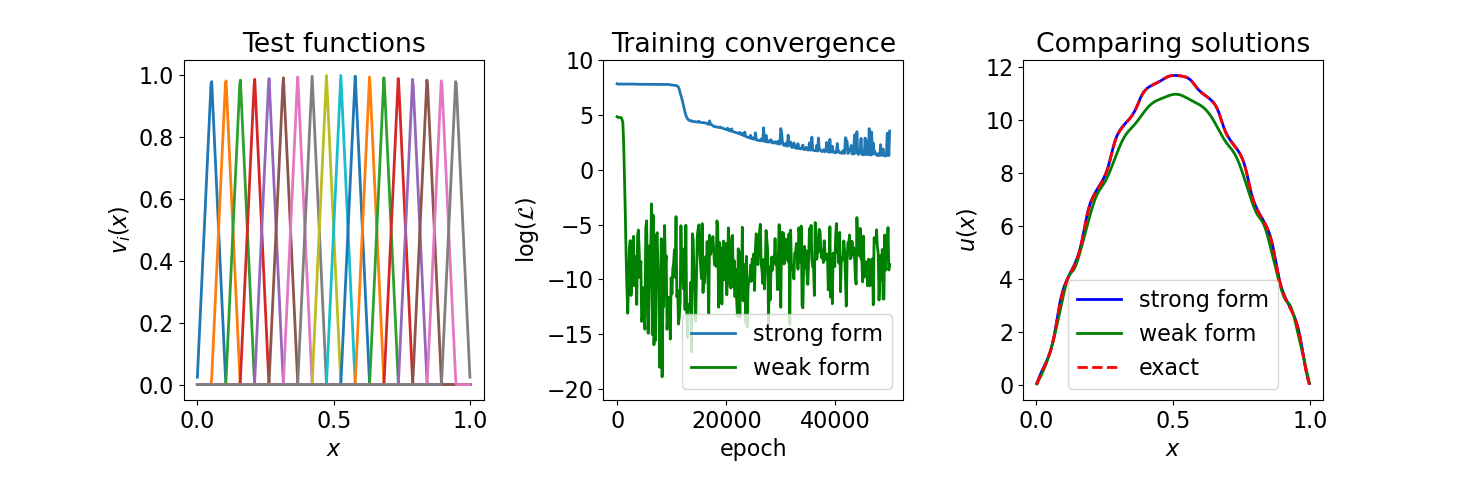}
\caption{The $N=18$ finite element hat functions are shown for reference (left), as well as the training dynamics of the strong and weak form losses (center) with the corresponding solutions (right). The Petrov-Galerkin weak form solutions exhibit significant error despite the expressivity of the neural network discretization.}
\label{weak_form}
\end{figure}

We discretize the solution with a two-hidden-layer neural network with each layer containing $50$ hidden units and use hyperbolic tangent activation functions. Both the strong and weak form losses are unconstrained and are minimized with ADAM optimization. The learning rate is $1 \times 10^{-3}$, the optimization is run for $5 \times 10^5$ epochs, and there are $N=18$ hat functions in the test set. See Figure \ref{weak_form} for the results from the strong and weak form solutions to the multiscale bar problem. The minimizer of the strong form loss reproduces the exact solution, whereas the zero of the weak form system leads to a solution with a spuriously large stiffness. Despite the expressivity of the neural network discretization, which is evidenced in the accuracy of the strong form solution, the weak form solution inherits the well-known failure of a coarse finite element solution to capture even the low-frequency component of the multiscale displacement \cite{leung_nh-pinn_2022}. The expressivity of the neural network is hindered by a poor choice of test functions. Thus, on top of the challenges of obtaining a suitable set of test functions for arbitrary two- and three-dimensional geometries, the weak form solution is only as good as the choice of test functions. We note that similar problems occur using spectral Fourier or Legendre polynomial test functions when the maximum frequency of the test function cannot resolve the multiscale component of the solution. We suggest that if the analyst has in hand a set of test functions that is i) adapted to the geometry and BCs of the problem and ii) sufficiently expressive to avoid issues such as those of Figure \ref{weak_form}, the neural network discretization offers few advantages, and a Bubnov-Galerkin method should be used with this set of test functions. Given the limitations of the energy objective and the reliance of the weak form on constructing a suitable test set, we believe the strong form formulation to be best-suited for a flexible PINN-based solution methodology. Of course, we do not doubt the ability of methods based on the weak form to accurately solve PDEs; we simply propose that the advantages of the neural network discretization are more salient in the case of the strong form. As such, we work with the strong form loss going forward.

\section{Comparison of methods for BC enforcement}
\label{sec:comparison}

\paragraph{} To motivate our discussion of BC enforcement, we introduce the following abstract scalar boundary value problem:
\begin{equation}\label{generix}
\begin{aligned}
    \mathcal{G}\Big( u(\mathbf{x}) \Big) + f(\mathbf{x}) = 0, \quad \mathbf{x} \in \omega \\
    u(\mathbf{x}) = g(\mathbf{x}), \quad \mathbf{x} \in \partial \omega_D, \\
    \nabla u \cdot \mathbf{\hat n} = t(\mathbf{x}), \quad \mathbf{x} \in \partial \omega_N,
\end{aligned}
\end{equation}
where $u(\mathbf{x})$ is the solution field, $\mathcal{G}(\cdot)$ is a differential operator, $\omega$ is the computational domain, $\partial \omega_D$ is the Dirichlet region of the boundary, $\partial \omega_N$ is the Neumann region, and $\mathbf {\hat n}$ is the normal vector. Although Robin conditions will appear in our numerical examples, we restrict attention here to Dirichlet and Neumann boundaries for clarity. We note that, despite the challenges in enforcing BCs with PINNs, a comprehensive survey and comparison of techniques available for Dirichlet and Neumann BC enforcement is lacking in the PINNs literature. An initial step in this direction is \cite{berrone_enforcing_2023}, where standard penalty methods for Dirichlet BCs are compared against distance functions and Nitsche's method using the weak formulation of the physics loss. Because we choose to work with the strong form loss, Nitsche's method---which is only used for Dirichlet boundaries---is ruled out as a candidate for BC enforcement. Our task is therefore to conduct a comprehensive survey of BC enforcement techniques for PINNs and to carry out a side-by-side comparison of those most promising for three-dimensional geometries. The results of this study will be used to establish our suggestions for best practices. 

\begin{figure}[hbt!]
\centering
\includegraphics[width=0.7\textwidth]{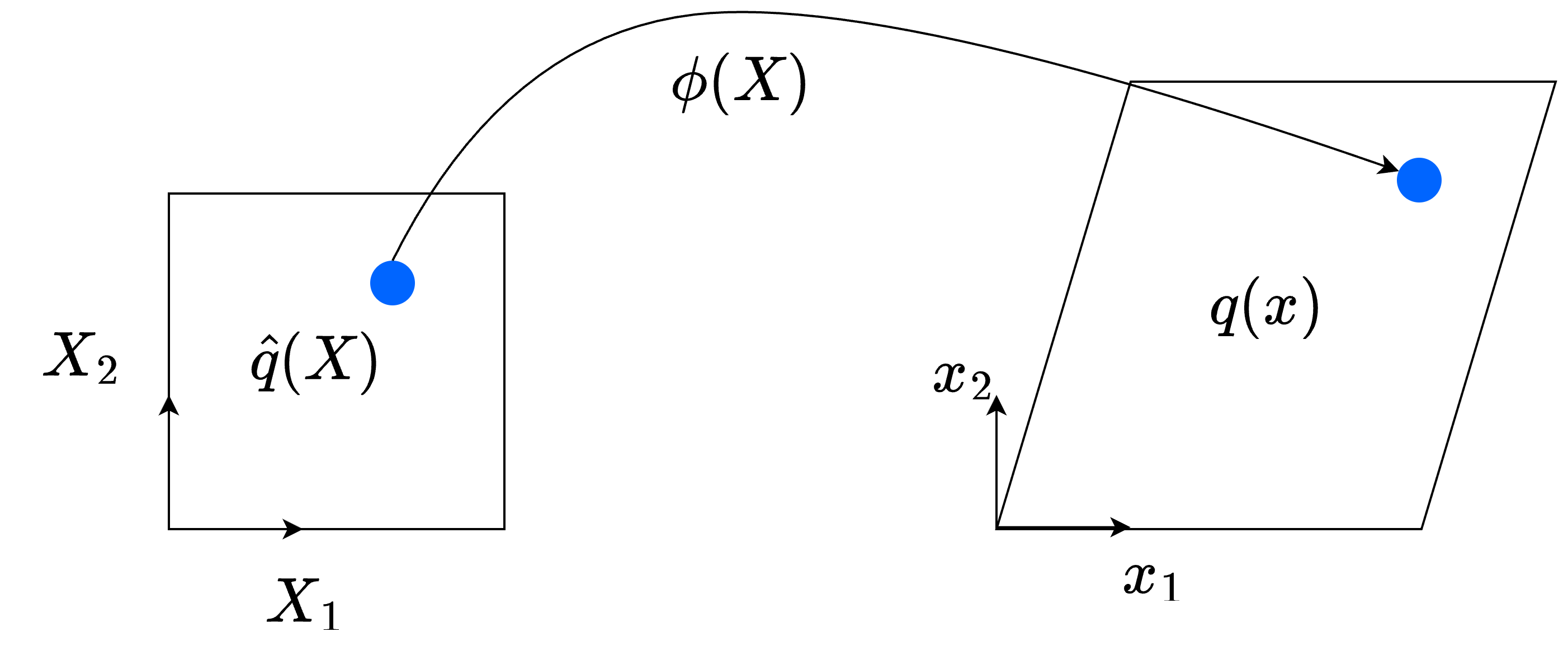}
\caption{The deformation $\boldsymbol \phi$ relates points $\mathbf{x}$ in the current geometry to corresponding points $\mathbf{X}$ in the reference geometry.}
\label{map}
\end{figure}

\paragraph{Coordinate transformation} We begin our survey with a discussion of coordinate transformation, which is a technique to map integrals and derivatives defined on one domain to their counterparts on another, often simpler, domain. We briefly describe coordinate transformation as a technique for BC enforcement. This is a topic which is best known in the scientific machine learning literature through \cite{li_fourier_2024}, in which coordinate transformation is used to map arbitrary geometries to tori in order to use the fast Fourier transform to compute convolution integrals. Let us assume that a quantity of interest $q(\mathbf{x})$ is defined over a complex ``current'' geometry given by $\mathbf{x} \in \omega$. Assume also that a simpler ``reference'' geometry is given by $\mathbf{X} \in \Omega$. The deformation map $\boldsymbol \phi: \mathbb{R}^3 \rightarrow \mathbb{R}^3$ relating the two geometries is given as $ \mathbf{x} = \boldsymbol \phi(\mathbf{X})$. See Figure \ref{map} for a schematic of the relation between the two configurations. Using the deformation, the change of variables formula for multivariate integration is

\begin{equation}\label{integral}
    \int_{\omega} q(\mathbf{x}) d\mathbf{x} = \int_{\Omega} q\Big( \boldsymbol{\phi}(\mathbf{X}) \Big) \text{det}\qty( \pd{\boldsymbol \phi}{\mathbf{X}}) d\mathbf{X} = \int_{\Omega} \hat q(\mathbf{X}) J(\mathbf{X}) d\mathbf{X},
\end{equation}

\noindent where $J(\mathbf{X})$ is the ``Jacobian determinant'' relating volume elements in the two configurations. We note that in the context of PINNs, the quantity of interest $q(\mathbf{x})$ will depend on spatial gradients with respect to the current coordinates. In this case, the derivatives also need to be transformed to the corresponding derivatives in the reference geometry. With this in mind, Eq. \eqref{integral} is a recipe for computing integrals of a quantity of interest defined in the current geometry entirely in terms of the reference geometry and a new field $\hat q(\mathbf{X})$ defined therein. For example, the strong form loss defined on the current geometry $\omega$ can be ``pulled back'' to a simpler geometry $\Omega$ to ease the enforcement of BCs. The difficulty with this coordinate transformation approach is threefold: first, we must obtain a deformation map $\boldsymbol \phi$ whose inverse transforms the computational domain $\omega$ back to a simple reference geometry (square/cube, circle/sphere, etc.). One way to do this is to use a neural network to represent the deformation and train it by solving another set of partial differential equations \cite{wang_pinn-mg_2025}. This additional PDE solve introduces significant overhead into the forward solve of interest, in addition to requiring the analyst to specify an appropriate reference geometry. Second, poor choices of reference geometry may lead to Jacobian determinant fields that, through their spatial variation, make the physics of the problem on the reference geometry more complicated. Third, though the reference geometry may be simpler, it is still necessary to devise a method to handle BCs there. We seek a technique for BC enforcement that is general enough to apply to the original problem geometry, and we thus eliminate coordinate transformation as a viable approach.

\begin{figure}[hbt!]
\centering
\includegraphics[width=0.7\textwidth]{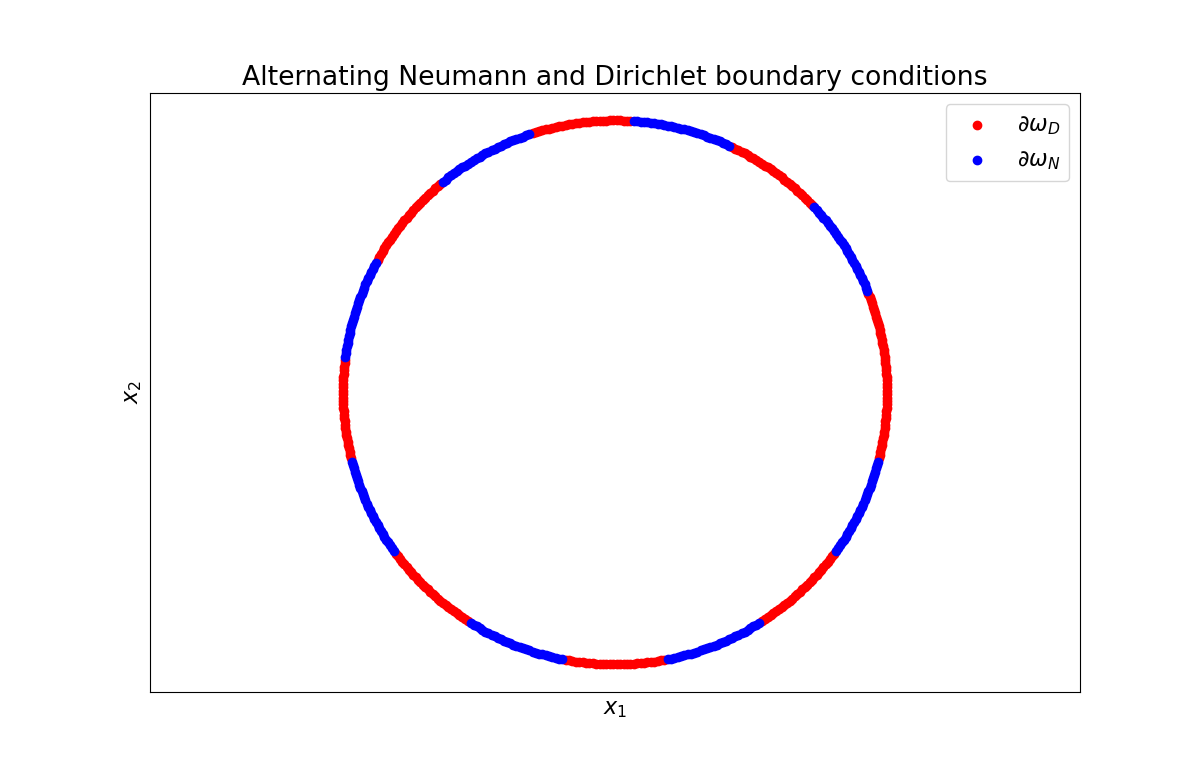}
\caption{Mixed boundaries present challenges for methods based on distance functions, as the distance function must be zero along $\partial \omega_D$ and non-zero along $\partial \omega_N$.}
\label{combo}
\end{figure}

\paragraph{Distance functions} Distance functions have become a popular technique for enforcing Dirichlet BCs with PINNs. Following \cite{wang_exact_2023, sukumar_exact_2022, sheng_pfnn_2021, berrone_enforcing_2023, luong_automatically_2024}, a distance function method discretizes the PDE solution as
\begin{equation}\label{distance}
    \hat u(\mathbf{x}; \boldsymbol \theta) = G(\mathbf{x}) + D(\mathbf{x}) \mathcal{N}(\mathbf{x} ; \boldsymbol \theta),
\end{equation}
where $\mathcal{N}(\mathbf{x};\boldsymbol \theta)$ is a neural network that need not satisfy the BCs, and the functions $D(\mathbf{x})$ and $G(\mathbf{x})$, respectively, satisfy $D(\mathbf{x}) = 0$ for $\mathbf{x} \in \partial \omega_D$ and $G(\mathbf{x}) = g(\mathbf{x})$ for $\mathbf{x} \in \partial \omega_D$. We note that the function $D(\mathbf x)$ does not need to be a true signed distance function, in that we do not require that $\lVert \nabla D \rVert=1$, which is the Eikonal property. The above construction ensures that the inhomogeneous Dirichlet BCs are enforced by construction. Note that as it stands, the distance function approach given in Eq. \eqref{distance} does not have any bearing on the Neumann BCs, which need to be handled separately. Furthermore, constructing the two functions $G(\mathbf{x})$ and $D(\mathbf{x})$ may be a significant bottleneck in the solution process. Consider the circular geometry shown in Figure \ref{combo}. The boundary is built up from alternating patches of Dirichlet and Neumann BCs. Analytic distance functions are not readily available for problems of this type, and standard methods used to numerically construct distance functions would need to be modified, given the interspersing of the Neumann patches \cite{sukumar_exact_2022}. As is shown in \cite{sheng_pfnn_2021}, separate neural networks can be trained to obtain $G(\mathbf{x})$ and $D(\mathbf{x})$, but this introduces additional overhead to the PDE solve and ambiguity as to what conditions these two functions should satisfy along the Neumann boundary. To remedy this ambiguity, a scalar distance function (e.g., a distance function that satisfies the Eikonal equation) could be built, which computes the distance to the nearest point on the nearest Dirichlet patch, thus ensuring that the function $D(\mathbf{x})$ is non-zero along the Neumann patches. However, having a disconnected set of boundary curves increases the complexity of obtaining the distance function and still leaves Neumann boundaries to be enforced separately. 

\paragraph{} While this method has not been explored as extensively as that of Eq. \eqref{distance} in the PINNs literature, the authors in \cite{sukumar_exact_2022} show that Neumann and Robin boundary conditions can also be enforced explicitly in the discretization of the PDE solution. For example, if the entire boundary is Neumann with prescribed normal gradient $q(\mathbf x)$, a discretization which satisfies this BC by construction is 
\begin{equation*}
    \hat u(\mathbf x; \boldsymbol \theta_1,\boldsymbol \theta_2) = \mathcal N_1(\mathbf x; \boldsymbol \theta_1) - D(\mathbf x) \nabla D(\mathbf x) \cdot \nabla \mathcal N_1(\mathbf x ; \boldsymbol \theta_1) - D(\mathbf x) q(\mathbf x) + D(\mathbf x)^2 \mathcal N_2(\mathbf x; \boldsymbol \theta_2),
\end{equation*}
\noindent where $\mathcal N_1$ and $\mathcal N_2$ are independent neural networks with parameters $\boldsymbol \theta_1$ and $\boldsymbol \theta_2$. Note that it is now necessary for $D(\mathbf x)$ to satisfy the Eikonal property, at least near the boundary of the domain. Though the discretization shown here assumes the entire boundary is Neumann, the authors modify this basic construction to account for domains with mixed BCs. However, this extension to mixed boundaries significantly increases in the complexity of the implementation. That being said, there are more fundamental challenges introduced by the strong enforcement of Neumann boundaries. First, we point out that there are now two networks that must be trained to obtain the PDE solution, which will increase the number of degrees of freedom in the optimization problem, and thus the cost. Second, the gradient of $\mathcal N_1$ is required in the discretization, which ensures that, given our reliance on the strong form loss, this network will be differentiated three times with respect to its inputs, further increasing the computational cost of the method. Third, offline calculations are required to extend the Neumann boundary $q(\mathbf x)$ into the domain, and to compute a distance function $D(\mathbf x)$ satisfying the Eikonal property near the boundary. We note that the requirement of offline calculations is not inherently a disadvantage, however in our opinion, the ability to quickly (in time) and efficiently (in lines of code) implement PINNs solutions is one of the appealing features of neural network discretizations. Thus, we consider the requirement of costly and/or complicated offline calculations as a disadvantage, though this is considered in the context of other pros and cons of each method. The authors provide methods to build distance functions of this sort from combinations of line segments, thus enabling the construction of $D(\mathbf x)$ for polygonal domains, but discuss to a lesser extent the treatment of domains with curved sides for which analytical distance functions are not known. A method is presented to handle arbitrary curved boundaries, but it is limited to domains in two spatial dimensions and further increases the complexity of the implementation \cite{sukumar_exact_2022}. To the best of our knowledge, this method has yet to be extended to three spatial dimensions, but we speculate that a parameterization of the boundary surface would be required if so, giving rise to additional offline steps to handle the BCs. The fourth and final challenge of hard enforcement of Neumann boundaries is that Eikonal distance functions have singular Laplacians when the boundary curve is not continuously differentiable. As the authors discuss, when the Eikonal property is not satisfied globally, these singularities are limited to points near the boundary. To remedy this issue, they avoid forming the physics-informed loss near the boundary with a careful choice of the integration grid, as the PINN solution is inaccurate if the singular region of the distance function is included in the physics-informed loss. We remark that a recent work modifies these methods to remove the singularity of the distance function for boundaries that are not continuously differentiable \cite{goschel_enforcing_2025}. In the authors' own words, this new approach ``[increases] the complexity of the implementation compared to weakly enforced boundary conditions \ldots [and] inference times [by] up to 30\%.'' In one numerical example, instead of a single trainable function for the discretized solution, there are $10$. We note that the cost of the optimization problem is not necessarily linear in the number of trainable functions, as the output layer of the network can be modified to be higher-dimensional, thus allowing each output function to share parameters of hidden layers. That being said, a higher-dimensional network output does bring with it additional parameters, as well as increased complexity of the optimization problem. Furthermore, this method has not been extended to three spatial dimensions.

\paragraph{} The above discussion shows that hard enforcement strategies for Neumann boundaries require offline processing to obtain Eikonal distance functions, while increasing the complexity of the implementation, the order of differentiation on the solution, and the number of trainable parameters. Perhaps most importantly, the Laplacian of the Eikonal distance function is singular for any domain with corners---an extremely common feature of engineering designs---leading to the failure of physics-informed training. While enforcing Dirichlet boundaries exactly with a discretization like Eq. \eqref{distance} is appealing, our analysis of Neumann boundaries leads us to conclude that distance functions are not a preferred strategy for a general-purpose, easy-to-use PINN solver, which should allow for arbitrary combinations of Dirichlet, Neumann, and Robin boundaries on complex three-dimensional geometries to be handled straightforwardly. The interested reader is referred to Appendix A for numerical examples supporting these points and further discussion.

\paragraph{PINN-FEM} A final strategy for exactly enforcing Dirichlet boundaries is PINN-FEM, which is a hybrid of a finite element discretization and PINNs relying on domain decomposition \cite{sobh_pinn-fem_2025}. As in the distance function methods, PINN-FEM relies on pre-processing---in this case, meshing---and only handles Dirichlet BCs. Furthermore, we stress again that relying on a mesh undermines one of the main advantages of PINNs, which is its meshless character. Accordingly, this method is also ruled out.

\begin{table}[hbt!]
    \centering
    \begin{tabular}{|c|c|c|}
        \hline
        \textbf{Method} & \textbf{Avoids pre-processing?} & \textbf{Neumann BC?} \\ \hline
        Coordinate transformation \cite{burbulla_physics-informed_2023} & \xmark & \xmark \\ \hline
        Distance function \cite{sukumar_exact_2022} & \xmark & \xmark \\ \hline
        PINN-FEM \cite{sobh_pinn-fem_2025} & \xmark & \xmark \\ \hline
        Standard penalty \cite{raissi_physics-informed_2019} & \cmark & \cmark \\ \hline
        Learning rate annealing \cite{wang_understanding_2020} & \cmark & \cmark \\ \hline
        SA-PINNs \cite{mcclenny_self-adaptive_2023} & \cmark & \cmark \\ \hline
        Lagrange multipliers \cite{makridakis_deep_2024} & \cmark & \cmark \\ \hline
        Augmented Lagrangian \cite{son_enhanced_2023} & \cmark & \cmark \\ \hline
    \end{tabular}
    \caption{A list of methods available to enforce BCs for PINN solutions based on the strong form loss.}
    \label{tab:methods}
\end{table}

 
\paragraph{} We have disqualified Nitsche's method because it does not apply to PINN solutions based on the strong form loss. We have also chosen to rule out methods that handle only Dirichlet boundaries and/or rely on significant pre-processing efforts to prepare the PDE solve. See Table \ref{tab:methods} for a list of BC enforcement methods for PINN formulations based on the strong form loss. Thus far, we have discussed the first three methods and disqualified them as approaches to a general-purpose PINN solver. We reiterate that this disqualification does not call into question the utility of these methods in specific contexts. Rather, we claim that these methods require pre-processing and/or do not apply to arbitrary combinations of Dirichlet and Neumann boundaries on three-dimensional geometries. The remaining five methods apply to both Dirichlet and Neumann boundaries without an offline processing stage: standard penalty methods \cite{raissi_physics-informed_2019}, learning rate annealing \cite{wang_understanding_2020}, SA-PINNs \cite{mcclenny_self-adaptive_2023}, Lagrange multipliers \cite{makridakis_deep_2024}, and the Augmented Lagrangian \cite{son_enhanced_2023}. In order to further explore the pros and cons of each of these methods, we perform a comparison on a two-dimensional test problem with inhomogeneous Dirichlet BCs. For our test problem, we solve the governing equation of forced steady-state heat conduction on the unit circle, given by
\begin{equation}\label{bvp}
\begin{aligned}
    \nabla^2 u(\mathbf{x}) + f(\mathbf{x}) = 0, \quad \mathbf{x} \in \omega,\\
    u(\mathbf{x}) = g(\mathbf{x}), \quad \mathbf{x} \in \partial \omega,
\end{aligned}
\end{equation}
where the entire boundary is Dirichlet such that $\partial \omega_D = \partial \omega$. For each method, we first establish the objective functions and the corresponding optimization problem. The temperature field is discretized by an MLP neural network $\mathcal{N}$ as $\hat u(\mathbf{x} ; \boldsymbol \theta) = \mathcal{N}(\mathbf{x} ; \boldsymbol \theta)$.

\paragraph{Standard penalty} A solution to Eq. \eqref{bvp} is obtained using a standard penalty method with
\begin{equation}\label{Lpenalty}
\begin{aligned}
\mathcal{L}^{\text{PENALTY}} (\boldsymbol \theta) = \frac{1}{2}\int_{\omega} \Big( \nabla^2 \hat u(\mathbf{x}) + f(\mathbf{x})\Big)^2 d\mathbf{x} + \frac{\beta}{2}\int_{\partial \omega} \Big(  \hat u(\mathbf{x}) - g(\mathbf{x})\Big)^2 ds,\\
\underset{\boldsymbol \theta}{\text{argmin }} \mathcal{L}^{\text{PENALTY}}(\boldsymbol \theta),
\end{aligned}
\end{equation}
\noindent where $\beta$ is a penalty hyperparameter \cite{raissi_physics-informed_2019}. This was the original approach for enforcing BCs with PINNs, and it remains a popular strategy owing to its simplicity.

\paragraph{Learning rate annealing} Depending on the choice of $\beta$, the two loss terms in Eq. \eqref{Lpenalty} may be unbalanced, leading to stiff gradient flow dynamics. To mitigate this issue, the penalty parameter can be dynamically updated with ``learning rate annealing'' \cite{wang_understanding_2020} per
\begin{equation}\label{lra}
\begin{aligned}
\mathcal{L}^{\text{LRA}} (\boldsymbol \theta) = \frac{1}{2}\int_{\omega} \Big( \nabla^2 \hat u(\mathbf{x}) + f(\mathbf{x})\Big)^2 d\mathbf{x} + \frac{\beta(t)}{2}\int_{\partial \omega} \Big(  \hat u(\mathbf{x}) - g(\mathbf{x})\Big)^2 ds,\\
\hat \beta_t = \frac{\text{max}\qty( \Big| \frac{1}{2}\nabla_{\boldsymbol \theta} \int_{\omega} \Big( \nabla^2 \hat u(\mathbf{x}) + f(\mathbf{x})\Big)^2 d\mathbf{x} \Big| )}{ \text{mean}\qty( \Big |\frac{1}{2} \nabla_{\boldsymbol \theta}\int_{\partial \omega} \Big(  \hat u(\mathbf{x}) - g(\mathbf{x})\Big)^2 ds \Big|)},\\
 \beta_t = (1-\alpha)\beta_t + \alpha \hat \beta_t,\\
\boldsymbol \theta_{t+1} = \boldsymbol \theta_t - \eta \nabla_{\boldsymbol \theta}\qty(  \frac{1}{2}\int_{\omega} \Big( \nabla^2 \hat u(\mathbf{x}) + f(\mathbf{x})\Big)^2 d\mathbf{x} + \frac{\beta_t}{2}\int_{\partial \omega} \Big(  \hat u(\mathbf{x}) - g(\mathbf{x})\Big)^2 ds ),
\end{aligned}
\end{equation}
\noindent where $0<\alpha<1$ defines a moving average of the penalty parameters over optimization iterations,  $\eta$ is the learning rate, and $t$ indexes the optimization epoch. Though we show a standard gradient descent update in Eq. \eqref{lra}, other optimizers, such as ADAM, can be used as well. 

\paragraph{SA-PINNs} Another modification to the standard penalty approach is SA-PINNs, which seeks to avoid the dependence of the accuracy of the BC enforcement on the choice of the penalty parameter. The objective function and optimization problem for SA-PINNs are given by 
\begin{equation}\label{sa-pinn}
\begin{aligned}
\mathcal{L}^{\text{SA}} (\boldsymbol \theta, \beta) = \frac{1}{2}\int_{\omega} \Big( \nabla^2 \hat u(\mathbf{x}) + f(\mathbf{x})\Big)^2 d\mathbf{x} + \frac{1}{2}\int_{\partial \omega} \beta(\mathbf{x})\Big(  \hat u(\mathbf{x}) - g(\mathbf{x})\Big)^2 ds,\\
\underset{\beta}{\text{argmax }}\underset{\boldsymbol \theta}{\text{argmin }} \mathcal{L}^{\text{SA}}(\boldsymbol \theta, \beta),
\end{aligned}
\end{equation}
\noindent where $ds$ is a differential area element on the boundary and $\beta(\mathbf{x})$ is a penalty field defined over the boundary which is driven up according to the pointwise boundary error \cite{mcclenny_self-adaptive_2023}. In our implementation, as well as that of the original paper, $\beta(\mathbf{x})$ is discretized at integration points. 

\paragraph{Lagrange multipliers} Departing from penalty methods, another min-max formulation of the PINN solution is that of Lagrange multipliers, with the objective function and optimization problem defined by
\begin{equation}\label{lagrange}
\begin{aligned}
\mathcal{L}^{L} (\boldsymbol \theta, \lambda) = \frac{1}{2}\int_{\omega} \Big( \nabla^2 \hat u(\mathbf{x}) + f(\mathbf{x})\Big)^2 d\mathbf{x} + \int_{\partial \omega} \lambda(\mathbf{x})\Big(  \hat u(\mathbf{x}) - g(\mathbf{x})\Big) ds,\\
\underset{\lambda}{\text{argmax }}\underset{\boldsymbol \theta}{\text{argmin }} \mathcal{L}^{L}(\boldsymbol \theta, \lambda),
\end{aligned}
\end{equation}
\noindent where $\lambda(\mathbf{x})$ is a Lagrange multiplier field defined over the boundary \cite{makridakis_deep_2024}. As in SA-PINNs, the Lagrange multiplier is discretized at integration points in the numerical implementation. Unlike Eqs. \eqref{sa-pinn}, the solution parameters can be chosen to minimize the loss function by making the signed boundary error $\hat u(\mathbf{x}) - g(\mathbf{x})$ more negative. This suggests a more complex optimization problem than that of SA-PINNs, for which descent directions of the boundary loss with respect to the penalty parameter correspond to increased BC accuracy. 

\paragraph{Augmented Lagrangian} A popular hybrid of penalty methods and Lagrange multipliers is the Augmented Lagrangian method. A solution is obtained by solving an iterative sequence of optimization problems given by
\begin{equation}\label{AL}
\begin{aligned}
    \mathcal{L}^{\text{AL}}( \boldsymbol \theta, \lambda_t, \beta_t) = \frac{1}{2}\int_{\omega} \Big( \nabla^2 \hat u(\mathbf{x}) + f(\mathbf{x})\Big)^2 d\mathbf{x} + \int_{\partial \omega} \lambda_t(\mathbf{x})\Big(  \hat u(\mathbf{x}) - g(\mathbf{x})\Big) ds +\frac{1}{2} \int_{\partial \omega} \beta_t(\mathbf{x})\Big(  \hat u(\mathbf{x}) - g(\mathbf{x})\Big)^2 ds  ,\\
    \boldsymbol \theta_t = \underset{\boldsymbol \theta}{\text{argmin }} \mathcal{L}^{\text{AL}}( \boldsymbol \theta, \lambda_t,\beta_t),\\
    \lambda_t(\mathbf{x}) \leftarrow \lambda_t(\mathbf{x}) + \beta_t(\mathbf{x})\Big( \hat u(\mathbf{x}; \boldsymbol \theta_t) - g(\mathbf{x}) \Big), \\
    \beta_t(\mathbf{x}) \leftarrow \gamma \beta_t(\mathbf{x}),
\end{aligned}
\end{equation}
\noindent where $\gamma >1$ defines the monotonically increasing update to the penalty parameter and the iteration over $t$, which is the index of the ``inner loop'' optimization problem, continues until user-defined convergence criteria are met \cite{basir_physics_2022}. 

\paragraph{} We have now laid out the five candidate formulations of the objective function to handle the inhomogeneous Dirichlet boundary. We first discuss the solution of the optimization problems given in Eqs. \eqref{Lpenalty}-\eqref{AL} with first-order methods, such as ADAM and gradient descent. We then turn our attention to second-order Newton-type methods. Note that both these classes of optimization methods require the gradient of the objective function. There are gradient-free approaches to solving these optimization problems, e.g., nature-inspired algorithms \cite{sulaiman_theoretical_2018, waseem_investigation_2020}, but we do not consider these in this work. Returning to the test problem of Eq. \eqref{bvp}, we take the geometry $\omega$ of the test problem to be the unit circle centered at the origin. The source term $f(\mathbf{x})$ and the inhomogeneous Dirichlet boundary $g(\mathbf{x})$---which spans the entirety of the boundary curve---will be determined from an assumed solution $u(\mathbf{x})$ and Eq. \eqref{bvp} using the method of manufactured solutions. The temperature field we assume is 
\begin{equation*}
    u(\mathbf{x}) = 10\Big( x_1x_2 + \sin(\pi x_1) \sin(\pi x_2) ( 1 - x_1^2 - x_2^2) \Big).
\end{equation*}

\begin{figure}[hbt!]
\centering
\includegraphics[width=0.7\textwidth]{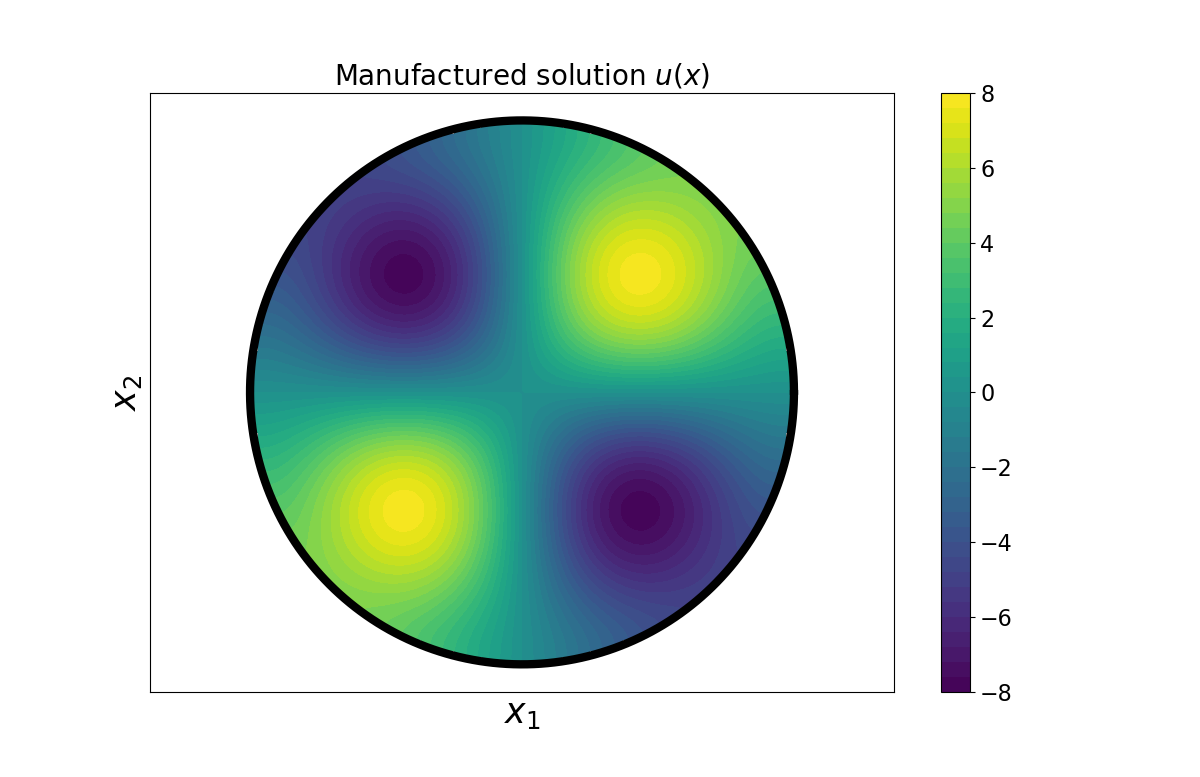}
\caption{The assumed temperature field used in the method of manufactured solutions, from which the source term and inhomogeneous Dirichlet BCs are computed.}
\label{temp_manufactured}
\end{figure}

See Figure \ref{temp_manufactured} to visualize this temperature field, which will be used throughout the study. We discretize the temperature field with a standard MLP neural network with two hidden layers, $25$ hidden units per layer, and hyperbolic tangent activation functions. We train the neural network corresponding to each objective function for $1 \times 10^4$ epochs using ADAM optimization for all gradient updates to the solution parameters $\boldsymbol \theta$. The learning rate for ADAM is set at $5 \times 10^{-3}$ for all methods except the Lagrange multiplier method of Eq. \eqref{lagrange}, for which it is set at $1\times10^{-3}$. A $100 \times 100$ background integration grid is set up on the square given by $[-1,1] \times [-1,1]$, and the points that lie inside the unit circle are taken as integration points. The boundary integration grid consists of $500$ equally spaced points around the circumference of the unit circle. We test the penalty method at two different penalty parameters: $\beta=1$ and $\beta=100$. For the learning rate annealing method, we set $\alpha=0.9$ per the suggestions of~\cite{wang_understanding_2020}. In the SA-PINNs method, the initial penalty parameters are all set to $1$, and the maximization is done with gradient descent with no momentum and a learning rate of $5 \times 10^{-1}$. Using gradient descent over ADAM ensures that the penalty stops increasing if the error is driven down to zero. In the min-max formulation of Lagrange multipliers, the initial values of the Lagrange multipliers are set to $0$, and the maximization is done with stochastic gradient descent without momentum with a learning rate of $1 \times 10^{-2}$. In the Augmented Lagrangian, we take the penalty update to be $\gamma=2$, initialize the penalty parameters at $1$, the Lagrange multipliers at $0$, and compute an update to the Lagrange multiplier and penalty parameter when the norm of the gradient of the objective falls within $1\%$ of its initial value. We note that these parameter settings are the consequence of experimentation to maximize the performance of each method.

\begin{figure}[hbt!]
\centering
\includegraphics[width=1.0\textwidth]{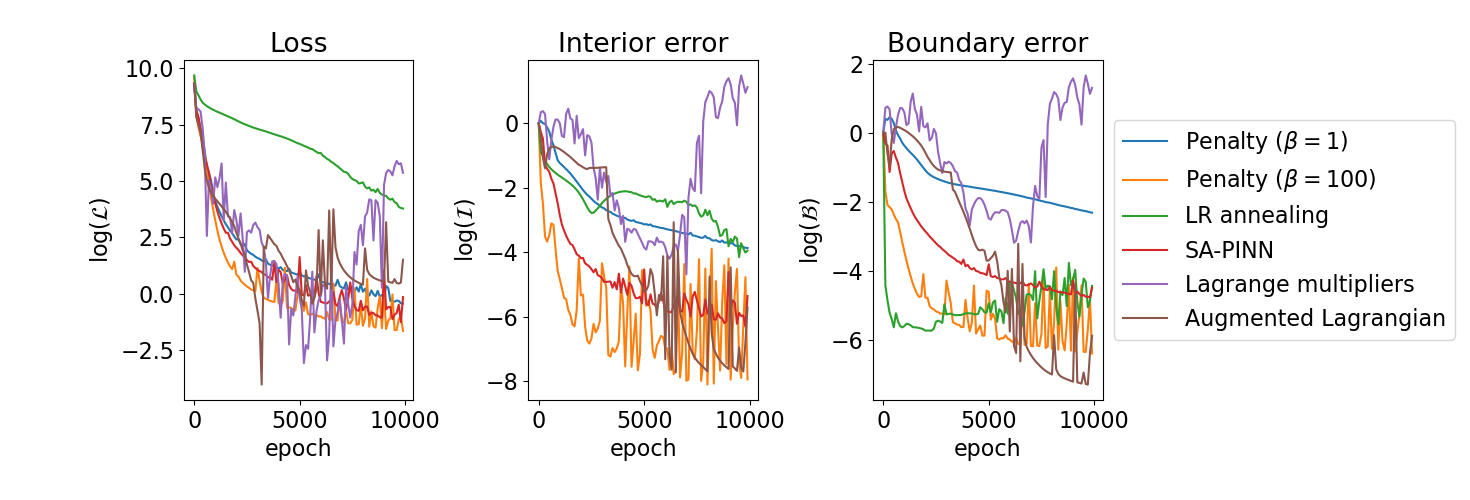}
\caption{Comparing the five formulations of the loss for enforcing inhomogeneous Dirichlet BCs, optimized with first-order methods. We report the absolute value of the loss with Lagrange multipliers, as it can take on negative values. Across the two error measures, the Augmented Lagrangian and SA-PINNs perform the best, with learning rate annealing also showing good performance.}
\label{study}
\end{figure}

\paragraph{} To compare the five formulations of the objective optimized with first-order methods, we define two different error measures, one for the interior error and one for the boundary error:
\begin{equation}\label{error}
    \mathcal{I}(\boldsymbol \theta) = \frac{\int_{\omega} | u(\mathbf{x}) - \hat u(\mathbf{x} ; \boldsymbol \theta)|d\mathbf{x}}{\int_{\omega} | u(\mathbf{x})|d\mathbf{x}}, \quad \mathcal{B}(\boldsymbol \theta) = \frac{\int_{\partial \omega} | g(\mathbf{x}) - \hat u(\mathbf{x} ; \boldsymbol \theta)|ds}{\int_{\partial \omega} | g(\mathbf{x})|ds}.
\end{equation}
Note that we choose the relative $L_1$ error as we believe it best captures an average pointwise difference between the two solution fields, which we interpret as a percent error. See Figure \ref{study} for the results of the comparison. When the penalty is set at a small value $(\beta=1)$, the standard penalty approach does not accurately enforce the BCs, and thus the interior error remains large. The larger penalty value ($\beta=100$) performs better, but it relies on a choice of hyperparameter which is specific to the problem at hand. For example, if the penalty is made too large, the optimizer will completely ignore the physics loss term. In contrast, by adaptively updating the magnitudes of the two loss terms, the learning rate annealing approach yields low boundary error without relying on a penalty hyperparameter. SA-PINNs performs well and relies only on the choice of the initial penalty parameter and the learning rate for the penalty parameter updates, which, along with the point-wise boundary error, determine the rate at which the optimizer drives them up. The min-max formulation of the Lagrange multiplier problem does not converge within the allotted training period, exhibiting oscillations on both long and short time scales. Similar behavior is observed using ADAM optimization for the Lagrange multiplier, as well as for any choice of learning rates we experimented with. Finally, the Augmented Lagrangian method with the given hyperparameters performs well on our test problem. In our test case, learning rate annealing, SA-PINNs, and the Augmented Lagrangian are the top performers. See Table \ref{tab:errors} for a comparison of the average interior and boundary error taken over the last $100$ steps of the optimization.

\begin{table}[h]
    \centering
    \begin{tabular}{|c|c|c|}
        \hline
        \textbf{Method} & \textbf{Interior error $\mathcal{I}$} & \textbf{Boundary error $\mathcal{B}$} \\
        \hline
        Penalty $(\beta=1)$ & $2.2 \times 10^{-2}$ & $9.8 \times 10^{-2}$ \\
        \hline
       Penalty $(\beta=100)$& $4.8 \times 10^{-3}$ & $5.6 \times 10^{-3}$ \\
        \hline
        LR annealing & $1.8 \times 10^{-2}$ & $9.0 \times 10^{-3}$\\
        \hline
        SA-PINN & $2.9 \times 10^{-3}$ & $8.7 \times 10^{-3}$ \\
        \hline
        Lagrange multipliers & $3.2 \times 10^{0}$ & $3.8 \times 10^{0}$ \\
        \hline
        Augmented Lagrangian & $1.2 \times 10^{-3}$& $1.4 \times 10^{-3}$ \\
        \hline
    \end{tabular}
    \caption{Comparing the techniques for BC enforcement by looking at the interior and boundary errors with the exact solution averaged over the last $100$ steps. The three best performers are learning rate annealing, SA-PINNs, and the Augmented Lagrangian method.}
    \label{tab:errors}
\end{table}

\begin{figure}[hbt!]
\centering
\includegraphics[width=1.0\textwidth]{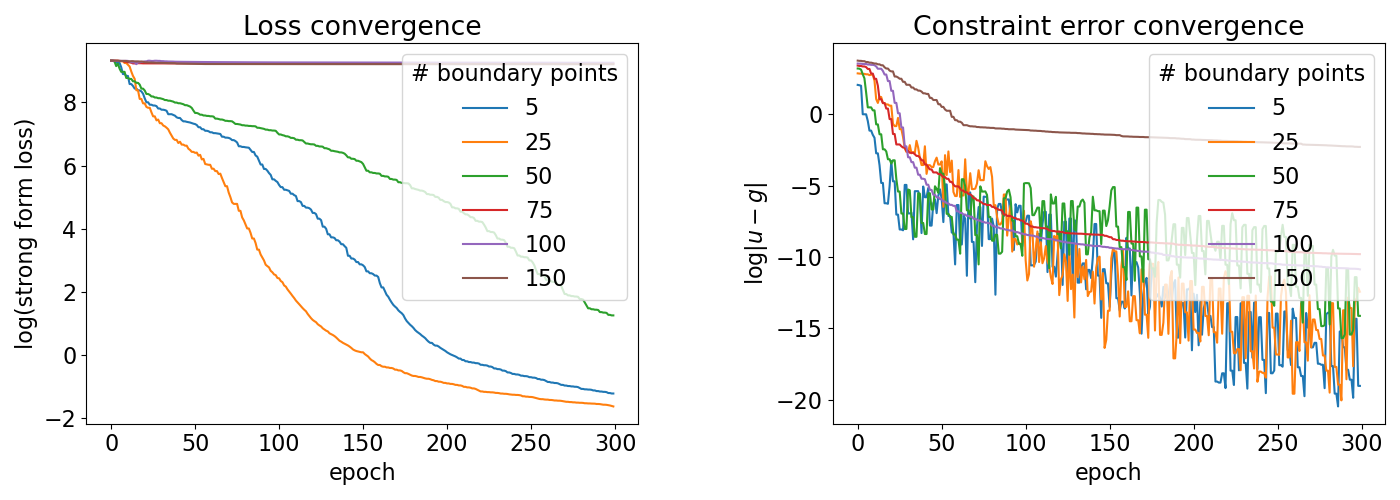}
\caption{Quasi-Newton constrained optimization with SciPy fails to converge the strong form loss with even a modest number of constraints.}
\label{conv}
\end{figure}

\paragraph{} We now discuss using second-order Newton-type methods to optimize the objective function. In particular, in place of the iterative min-max approach to the Lagrange objective of Eq. \eqref{lagrange}, we explicitly solve the nonlinear system corresponding to stationarity. This avoids the problematic min-max formulation at the cost of requiring second derivative information to perform updates to the solution parameters and Lagrange multiplier. The objective function is the Lagrange objective of Eq. \eqref{lagrange}, but we denote it with superscript ``$N$'' to distinguish this different approach to finding a stationary point. Accordingly, the stationarity condition and the Newton updating are written as
\begin{equation}\label{newton}
\begin{aligned}
\mathcal{L}^{\text{N}} (\boldsymbol \theta, \lambda) = \frac{1}{2}\int_{\omega} \Big( \nabla^2 \hat u(\mathbf{x}) + f(\mathbf{x})\Big)^2 d\mathbf{x} + \int_{\partial \omega} \lambda(\mathbf{x})\Big(  \hat u(\mathbf{x}) - g(\mathbf{x})\Big) ds,\\ \pd{ \mathcal{L}^{\text{N}} }{\boldsymbol \lambda} = \mathbf{0}, \quad \pd{\mathcal{L}^{\text{N}}}{\boldsymbol \theta} = \mathbf{0}, \\
\begin{bmatrix}
    \boldsymbol \theta \\ \boldsymbol \lambda
\end{bmatrix}^{t} \leftarrow \begin{bmatrix}
    \boldsymbol \theta \\ \boldsymbol \lambda
\end{bmatrix}^{t} - \begin{bmatrix} \partial^2 \mathcal{L}^{\text{N}}/ \partial \boldsymbol \theta \partial  \boldsymbol \theta & \partial^2 \mathcal{L}^{\text{N}}/ \partial \boldsymbol \theta \partial \boldsymbol \lambda \\ \partial^2 \mathcal{L}^{\text{N}}/ \partial \boldsymbol \lambda \partial \boldsymbol \theta & \partial^2 \mathcal{L}^{\text{N}}/ \partial \boldsymbol \lambda \partial \boldsymbol \lambda
\end{bmatrix}^{-1} \begin{bmatrix} \partial \mathcal{L}^{\text{N}} / \partial \boldsymbol \theta \\ \partial \mathcal{L}^{\text{N}} / \partial \boldsymbol \lambda
\end{bmatrix},
\end{aligned}
\end{equation}
\noindent where $t$ indexes the Newton step, and with slight abuse of notation, we keep the Lagrange multiplier continuous in the objective function but use $\boldsymbol \lambda$ (bold) to indicate its discretization at boundary integration points. The Newton updating is continued until the norm of the residual of the system for stationarity falls below a user-specified threshold. As we discuss shortly, an approximation of the Hessian matrix is often used in practice. Second-order methods have received more attention in the context of neural network training in recent years. In \cite{korbit_exact_2024}, Newton optimization is shown to accelerate the training of deep neural networks. In the context of PINNs, second-order optimizers have been systematically compared and shown to provide accelerated convergence \cite{kiyani_optimizing_2025, urban_unveiling_2025}. One known issue with Newton methods is that the updates converge to general stationary points in the loss landscape, rather than explicitly looking for minima \cite{dauphin_identifying_2014}. Numerical experimentation solving PINN-based optimization problems has suggested that true Newton methods, where the Hessian matrix in Eq. \eqref{newton} is computed exactly, often exhibit such pathological behavior \cite{rowan_nonlinear_2025}. As such, we choose to avoid true Newton methods in favor of quasi-Newton methods, where the Hessian is approximated using the BFGS methodology \cite{nocedal_quasi-newton_1999}. Returning to the Newton formulation of the test problem given in Eq. \eqref{newton}, we first experiment with the trust constrained SciPy optimizer, which uses the BFGS approximation of the Hessian to build a local quadratic model of the objective but limits the optimization step to a region whose size is adaptively updated based on the accuracy of the quadratic model. We find that when the number of constraints (the number of boundary integration points, in other words) is very small, the optimizer converges to a solution for which the strong form loss is approximately zero and the BCs are accurately satisfied. However, when the number of boundary integration points takes on realistic values of $100$-$500$, the optimizer gets stuck at solutions with large values of the strong form loss. In Figure \ref{conv}, we sweep over different sizes of the boundary constraints and show that the optimizer reliably drives the constraint error down to zero, but that the strong form loss stagnates once the number of boundary integration points becomes too large. We note that the neural network discretization consists of two hidden layers, each with $30$ hidden units, corresponding to $1050$ total parameters. This ensures that the problem is not over-constrained, given that the number of constraints in the study is at most $150$. Similar results are observed for the Sequential Least Squares Quadratic Programming (SLSQP) optimizer, which iteratively approaches a stationary point of the Lagrange function using a quadratic approximation of the objective with the BFGS Hessian, a linear approximation of the constraints, and line search. This is the one other quasi-Newton method in SciPy that handles constraints.

\paragraph{Summary} We ruled out coordinate transformations, distance functions, and PINN-FEM on first principles. These methods do not handle Neumann boundary conditions and require significant pre-processing steps. We then studied the standard penalty approach, learning rate annealing, SA-PINNs, Lagrange multipliers, and the Augmented Lagrangian method. We found that Lagrange multipliers with a second-order quasi-Newton optimizer struggled to handle the large number of nonlinear constraints. The study of Figure \ref{study} suggests that of the five methods solved with first-order optimization, the three best performing are SA-PINNs, the Augmented Lagrangian method, and learning rate annealing. Even when run for up to 30,000 epochs, the min-max formulation of Lagrange multipliers does not converge to a solution under any hyperparameter settings we could find, including choice of first-order optimization strategy, and is thus disqualified as a technique for BC enforcement. While the penalty method with the larger penalty parameter also performs well, we believe SA-PINNs, the Augmented Lagrangian, and learning rate annealing offer certain advantages. All of them involve adaptive balancing of the loss terms to make expensive hyperparameter tuning unnecessary. Thus, on the basis of its sensitivity to hyperparameters, we choose to rule out the standard penalty method. However, in our experience, SA-PINNs can also exhibit sensitivity to hyperparameters, namely, to the learning rate for the penalty parameters. If this learning rate is too low, SA-PINNs require a substantial number of optimization steps for the penalty parameters to stabilize. This increased convergence time is noted in \cite{hu_conditionally_2025}. A second problem with SA-PINNs is the introduction of additional trainable parameters, which balloons the number of optimization variables in two and three spatial dimensions. While the authors in \cite{son_enhanced_2023} remark that the method introduces $10,000$ trainable parameters in two dimensions, we remark that at our chosen resolution level, SA-PINNs would introduce up to $50,000$ additional optimization variables. With the neural network architectures we use, treating BCs with SA-PINNs corresponds to a fifty-fold increase in the size of the optimization problem. For these reasons, we narrow the focus of subsequent numerical examples to learning rate annealing and the Augmented Lagrangian method. We choose to emphasize the Augmented Lagrangian given its superior performance in our study. However, each numerical example also provides a comparison to the learning rate annealing method. We note that both methods require storing quantities at boundary integration points. In the case of learning rate annealing, the single penalty field is stored and updated throughout optimization. With the Augmented Lagrangian method, both the penalty and Lagrange multiplier fields are stored on the surface. Neither method requires any off-line computations. In the following section, we lay out a strategy for integrating the Augmented Lagrangian method for BC enforcement with arbitrary three-dimensional geometries.


\section{Solution framework for PINNs on 3D geometries}
\label{sec:method}

\paragraph{} Our task is to illustrate a methodology for solving general linear or nonlinear boundary value problems with any combination of Dirichlet, Neumann, and Robin BCs on arbitrarily complex three-dimensional geometries using PINNs. This task is primarily one of synthesis---we have weighed the pros and cons of different formulations of the physics loss, and concluded that the strong form loss is the best-equipped for our goal of a general-purpose PINN-based solver. Similarly, though the BC enforcement strategies we rely on have been proposed in the literature previously, a comprehensive comparison has been lacking. Having committed to the strong form loss and having narrowed the set of viable constraint enforcement strategies with such a comparison, the remaining ingredient in the numerical solution is the geometry handling. Once again, we do not introduce any novel tools to this end. Rather, we show how existing strategies can be integrated into a PINN solution workflow. We opt to use level sets to generate three-dimensional geometries, as this offers a streamlined and customizable approach to geometry generation. We note, however, that this choice is made with no loss of generality. Any standard geometry file will come equipped with the information required to build the interior and surface integration grids we extract from the level set. In our construction, all level sets define geometries which are subsets of the unit cube given by $\omega^B = [0,1]^3$, where the superscript ``$B$'' refers to the ``background'' domain. Consider a scalar-valued level set function $\Phi(\mathbf{x})$. The bulk $\omega$ and the boundary surface $\partial \omega$ of the three-dimensional geometry are defined by 
\begin{equation*}
    \omega = \{ \mathbf{x} | \Phi(\mathbf{x}) < 0\},  \quad \partial\omega = \{ \mathbf{x} | \Phi(\mathbf{x}) = 0\}.
\end{equation*}

\begin{figure}[hbt!]
\centering
\includegraphics[width=0.95\textwidth]{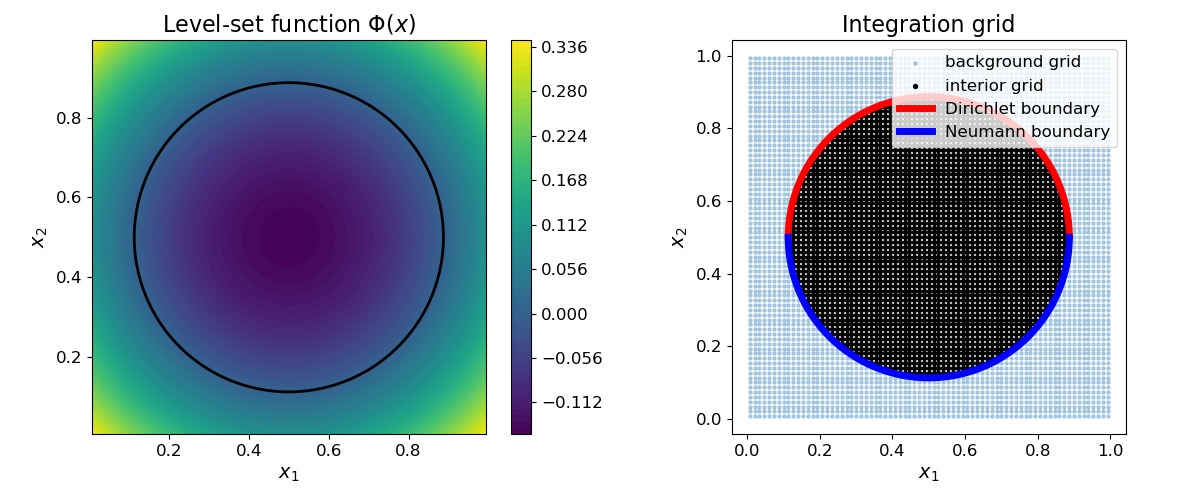}
\caption{As an example, a two-dimensional level set function is shown with its zero isocontour in black (left). The level set function is used to build an interior integration grid from a background integration, and the boundary is split up into Dirichlet and Neumann regions (right).}
\label{level_set_schematic}
\end{figure}

We define a uniform grid of points on the background domain as $\mathcal{P} = \{ \mathbf x_i \}_{i=1}^{N_B} \in \Omega^B$ where $N_B$ is the number of integration points. Next, we build an ``interior integration grid'' $\mathcal{X}$ using the level set with
\begin{equation*}
    \mathcal{X} = \{ \mathbf x \in \mathcal{P} | \Phi( \mathbf x) < 0\}.
\end{equation*}
This allows us to approximate integrals over the computational domain with
\begin{equation*}
    \int_{\omega} f(\mathbf{x}) d\mathbf{x} \approx \sum_{i=1}^N \Delta V f( \mathcal{X}_i),
\end{equation*}
\noindent where $N$ is the number of interior integration points, the set of interior integration points is $\mathcal{X}=\{ \mathcal{X}_i\}_{i=1}^N$, and $\Delta V=1/N$ is the volume element defined by the uniform integration grid. To enforce BCs, we must also compute surface integrals. This is accomplished with the marching cubes algorithm, which provides a triangulation of the zero isocontour of the level set evaluated on the background integration grid \cite{lorensen_marching_1987}. Returned from the marching cubes algorithm are sets of vertex points and the triangular faces formed by these vertex points. We define surface integration points as the centroids of the triangular faces, and compute outward-facing normal vectors by taking the cross product of the vectors defining the edges of the triangular faces. This allows us to approximate surface integrals as 
\begin{equation*}
    \int_{\partial \omega} \mathbf{q}(\mathbf{x}) \cdot \mathbf{ \hat n } ds \approx \sum_{j=1}^M \mathbf{q}( \mathcal{S}_j) \cdot \mathbf{\hat n}_j \Delta A_j,
\end{equation*}
\noindent where $M$ is the number of surface integration points, $\mathcal{S} = \{ \mathcal{S}_j \}_{j=1}^M$ is the surface integration grid obtained from marching cubes, $\mathbf{\hat n}_j$ is the outward facing unit normal corresponding to the $j$-th triangular facet, and $\Delta A_j$ is the corresponding area element. The particular implementation we use is Python's $\texttt{skimage}$ marching cubes algorithm. We remark that there are more sophisticated algorithms for integrating over volumes and surfaces defined by level sets \cite{saye_high-order_2015, hubrich_numerical_2019}. In this work, we rely on the grid of background points being sufficiently dense to obtain accurate volume integrals using a uniform weight $\Delta V$. The surface integration grid will be further broken down into the Dirichlet and Neumann components, which are given by $\mathcal{S}^D = \{ \mathcal{S}_j^D \}_{j=1}^{M_D}$ and $\mathcal{S}^N = \{ \mathcal{S}_k^N \}_{k=1}^{M_N}$ respectively. The superscripts ``$D$'' and ``$N$'' will be used similarly to index the area elements and normal vectors belonging to the two boundary sets. Figure \ref{level_set_schematic} illustrates a level set function and its zero isocontour, as well as background and interior integration grids, and a boundary split up into Dirichlet and Neumann regions. With the interior integration grid, the boundary integration grid, and the corresponding surface normals, we illustrate a PINN-based solution methodology for BCs on arbitrarily complex three-dimensional geometries using the Augmented Lagrangian method. We discretize the solution and leave the option open to enforce homogeneous Dirichlet BCs on a subset of the domain with a distance function $D(\mathbf{x})$. Thus, the discretization of the solution field is given by 

\begin{equation*}
    \hat u(\mathbf{x} ; \boldsymbol \theta) = D(\mathbf{x}) \mathcal{N}(\mathbf{x}; \boldsymbol \theta).
\end{equation*}

With the Lagrange multiplier field $\lambda(\mathbf{x})$ and penalty field $\beta(\mathbf{x})$ discretized at surface integration points, the discretized Augmented Lagrangian objective at iteration $t$ for the generic scalar BVP of Eq. \eqref{generix} is
\begin{equation}\label{ALdis}
\begin{aligned}
    \mathcal{\hat L}^{\text{AL}}( \boldsymbol \theta, \boldsymbol \lambda_t, \boldsymbol \beta_t) = \frac{1}{2} \sum_{i=1}^N \Delta V\Big( \mathcal{G}\Big( \hat u(\mathbf{ \mathcal{X}}_i)\Big) + f(\mathbf{\mathcal{X}}_i)\Big)^2  + \sum_{j=1}^{M_D} \Delta A^D_j \lambda^D_{j,t}\Big(  \hat u(\mathcal{S}^D_j) - g(\mathcal{S}^D_j)\Big)  \\+\frac{1}{2} \sum_{j=1}^{M_D}\Delta A^D_j\beta^D_{j,t}\Big(  \hat u(\mathcal{S}^D_j) - g(\mathcal{S}^D_j)\Big)^2 
    + \sum_{k=1}^{M_N}\Delta A_k^N \lambda^N_{k,t}\qty( \nabla \hat u (\mathcal{S}^N_k)\cdot \mathbf{\hat n}^N_k - t(\mathcal{S}^N_k ))+ \\ \frac{1}{2}\sum_{k=1}^{M_N}\Delta A_k^N \beta^N_{k,t}\qty(\nabla \hat u(\mathcal{S}^N_k) \cdot \mathbf{\hat n}^N_k - t(\mathcal{S}^N_k))^2.
\end{aligned}
\end{equation}


\begin{table}[h]
    \centering
    \begin{tabular}{|c|c|}
        \hline
        \textbf{Hyperparameter} & \textbf{Recommended setting} \\
        \hline
        Learning rate & $(1 \times 10^{-3})-(5 \times 10^{-3})$ \\
        \hline
        Initial penalty $\boldsymbol \beta_0$ & $\mathbf 1$ \\
        \hline
        Initial Lagrange multiplier $\boldsymbol \lambda_0$ & $\mathbf 0$ \\
        \hline
        Rate of penalty increase $\gamma$ & $2$ \\
        \hline
        Objective convergence $\mathcal{Z}_f$ & $(2.5 \times 10^{-3}) - (1\times10^{-2})$ \\
        \hline
        Boundary convergence $\mathcal{B}_f$ & $(5 \times 10^{-3}) - (1\times10^{-2})$  \\ \hline
        Gradient convergence $\mathcal{R}_f$ & $(1 \times 10^{-3}) - (1\times10^{-2})$ \\ \hline
    \end{tabular}
    \caption{Hyperparameters used for the PINN-based solution framework based on the Augmented Lagrangian method. The learning rate is used in ADAM optimization for the inner-loop problem to find a minimum of the Augmented Lagrangian objective at iteration $t$. The initial solution parameters $\boldsymbol \theta_0$ are obtained with the built-in Xavier initialization in PyTorch.}
    \label{tab:hyper}
\end{table}

\begin{algorithm}
\caption{Augmented Lagrangian for PINNs}\label{alg:AL}
    \begin{algorithmic}[1]
        \Require $\boldsymbol \theta_0, \boldsymbol \beta_0,\boldsymbol \lambda_0$ \Comment{initial parameters}
        \Require $\gamma, \mathcal{B}_f, \mathcal{Z}_f, 
        \mathcal{R}_f$ \Comment{penalty update and convergence criteria}
        \State $\mathcal{Z}_0 = \mathcal{L}^{\text{AL}}( \boldsymbol \theta_0 , \boldsymbol \lambda_0 , \boldsymbol \beta_0)$ \Comment{initial objective value}
        \State $\mathcal{R}_0 = \lVert \nabla_{\boldsymbol \theta} \mathcal{L}^{\text{AL}}( \boldsymbol \theta_0 , \boldsymbol \lambda_0 , \boldsymbol \beta_0) \rVert^2$ \Comment{initial gradient value}
        \While{ $\mathcal{B}_t > \mathcal{B}_f \textbf{ and } \mathcal{Z}_t/\mathcal{Z}_0 > \mathcal{Z}_f$ } \Comment{outer loop}
        \While{$\lVert \partial \mathcal{L}^{\text{AL}}(\boldsymbol \theta, \boldsymbol \lambda_t,\boldsymbol \beta_t) / \partial \boldsymbol \theta \rVert > \mathcal{R}_f/\mathcal{R}_0$} \Comment{inner loop}
        \State ADAM optimization on $\mathcal{L}^{\text{AL}}(\boldsymbol \theta, \boldsymbol \lambda_t,\boldsymbol \beta_t)$
        \EndWhile
        \State $ \boldsymbol \lambda^D_{t+1} = \boldsymbol \lambda^D_t  + \boldsymbol \beta^D_t \odot \Big( \hat u(\mathcal{S}^D; \boldsymbol \theta_t) - g(\mathcal{S}^D) \Big), \quad \boldsymbol \beta_{t+1}^D = \gamma \boldsymbol \beta^D_t $ \Comment{update for Dirichlet boundary}
        \State $ \boldsymbol \lambda^N_{t+1} = \boldsymbol \lambda^N_t  + \boldsymbol \beta^N_t \odot \Big( \nabla \hat u(\mathcal{S}^N; \boldsymbol \theta_t) \cdot \mathbf{\hat n}^N - t(\mathcal{S}^N) \Big), \quad \boldsymbol \beta_{t+1}^N = \gamma \boldsymbol \beta^N_t $ \Comment{update for Neumann boundary}
        \State $ \boldsymbol \lambda_{t+1} = [ \boldsymbol \lambda^D_{t+1} , \boldsymbol \lambda^N_{t+1}], \quad \boldsymbol \beta_{t+1} = [ \boldsymbol \beta^D_{t+1} , \boldsymbol \beta^N_{t+1}] $ 
        \State $\mathcal{Z}_{t+1} = \mathcal{L}^{\text{AL}}(\boldsymbol \theta_t,\boldsymbol \lambda_t, \boldsymbol \beta_t)$ \Comment{current objective value}
        \State $\mathcal{B}_{t+1} = \mathcal{B}(\boldsymbol \theta_t) $ \Comment{current boundary error}
        \State $t = t + 1$
        \EndWhile
    \end{algorithmic}
\end{algorithm}

We remark that the same strategy for computing surface integrals to form penalties on the BCs is used for the learning rate annealing method. We use ADAM optimization to solve the inner loop Augmented Lagrangian problem with no mini-batching of the integration grid. As is standard for the algorithm, the initial penalty parameter $(\boldsymbol \beta_0)$ is set at a small value, and the initial Lagrange multiplier ($\boldsymbol \lambda_0$ ) is set to zero. The rate of increase of the penalty parameter ($\gamma$) is set at $2$, which is also standard for the method. The most consequential hyperparameters for the Augmented Lagrangian method are the convergence criteria. We define convergence of the inner loop problem to occur when the gradient falls below $\mathcal{R}_f \%$ of its initial value. Convergence of the optimization concludes when this gradient condition is met, when the boundary loss falls below $\mathcal{B}_f$, and when the objective is below $\mathcal{Z}_f\%$ of its initial value. See Table \ref{tab:hyper} for the list of hyperparameters that specify this algorithm and Algorithm \ref{alg:AL} for the breakdown of our implementation of the Augmented Lagrangian method. Note that the notation $\odot$ indicates element-wise multiplication. Numerical experimentation on the forthcoming example problems has informed our choice of recommended settings given in Table \ref{tab:hyper}. In the next section, we showcase the efficacy of this method on three problems from engineering mechanics. We also compare its performance to learning rate annealing. All computations are performed on an AMD Ryzen 9 5950X 16-Core Processor and 64GB unified memory, running Python 3.11.5 and PyTorch 2.1.2.

\section{Numerical examples}
\label{sec:examples}

\subsection{Fisher-KPP equation on two-way branch}

\paragraph{} The first test case of the proposed method is a scalar-valued nonlinear static PDE. The Fisher-Kolmogorov-Petrovsky-Psikunov (Fisher-KPP) equation is a nonlinear reaction-diffusion PDE used to model population growth and traveling wave fronts \cite{rohrhofer_approximating_2025}. The dynamic Fisher-KPP equation is given by 
\begin{equation}\label{kppt}
    \pd{u}{t} = \mu \nabla^2 u(\mathbf{x}) + r u(\mathbf{x})(1 - u(\mathbf{x})) + f(\mathbf{x}),
\end{equation}
where $u$ is the concentration of a chemical substance or population density, $f(\mathbf{x})$ is a source term, $\mu$ is the diffusivity, and $r$ is a reaction rate coefficient. Going forward, we set $\mu=1, r=1/2$ and work with the steady state of Eq. \eqref{kppt} which is given by 
\begin{equation}\label{kpp}
\begin{aligned}
    \nabla^2 u(\mathbf{x}) + r u(\mathbf{x})(1 - u(\mathbf{x})) + f(\mathbf{x}) = 0, \quad \mathbf{x} \in \omega, \\
    u(\mathbf{x}) = g (\mathbf{x}), \quad \mathbf{x} \in \partial \omega.
\end{aligned}
\end{equation}

The geometry $\omega$ of the problem is that of a two-way branch, inspired by an arterial geometry from biomedical imaging. The following level set can be used to generate this geometry:
\begin{equation}\label{branch}
    \Phi(\mathbf{x}) = \Big( \cos 2\pi x_1 - (1-2x_3)\Big)^2 + 9(x_2-0.5)^2 - 0.5 .
\end{equation}

\begin{figure}[hbt!]
\centering
\includegraphics[width=0.6\textwidth]{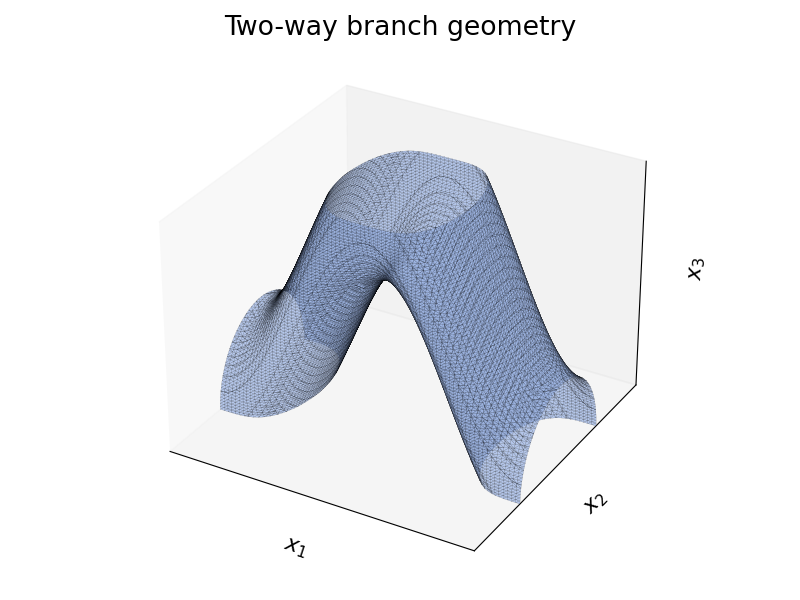}
\caption{Two-way branch geometry for the steady-state Fisher-KPP problem.}
\label{ex1_geo}
\end{figure}

The geometry is then defined by $\omega = \{ \mathbf{x} | \Phi(\mathbf{x}) < 0 \}$, as shown in Figure \ref{ex1_geo}. Recall that all geometries are subsets of the unit cube $[0,1]^3$. We use the marching cubes algorithm to mesh the surface defined by the zero isocontour of the level set given by Eq. \eqref{branch}. We enforce inhomogeneous Dirichlet BCs defined over the surface integration grid and a manufactured solution. The assumed solution is 
\begin{equation*}
    u(\mathbf{x}) = 10 \sin(3\pi x_3) \sin(2\pi x_1) \sin(\pi x_1) \sin(\pi x_2) \sin(\pi x_3),
\end{equation*}
which enforces that $u(\mathbf{x})=0$ along the surface of the unit cube in which the level set is defined. This is relevant, as the zero of the level set extends outside the unit cube; thus, the boundaries defined by the surface of the unit cube will not be picked up by marching cubes. These BCs will be handled in the discretization of the solution field per 
\begin{equation*}
    \hat u(\mathbf{x} ; \boldsymbol \theta) = \sin(\pi x_1) \sin(\pi x_2) \sin(\pi x_3)\mathcal{N}(\mathbf{x} ; \boldsymbol \theta).
\end{equation*}

\begin{figure}[hbt!]
\centering
\includegraphics[width=1.0\textwidth]{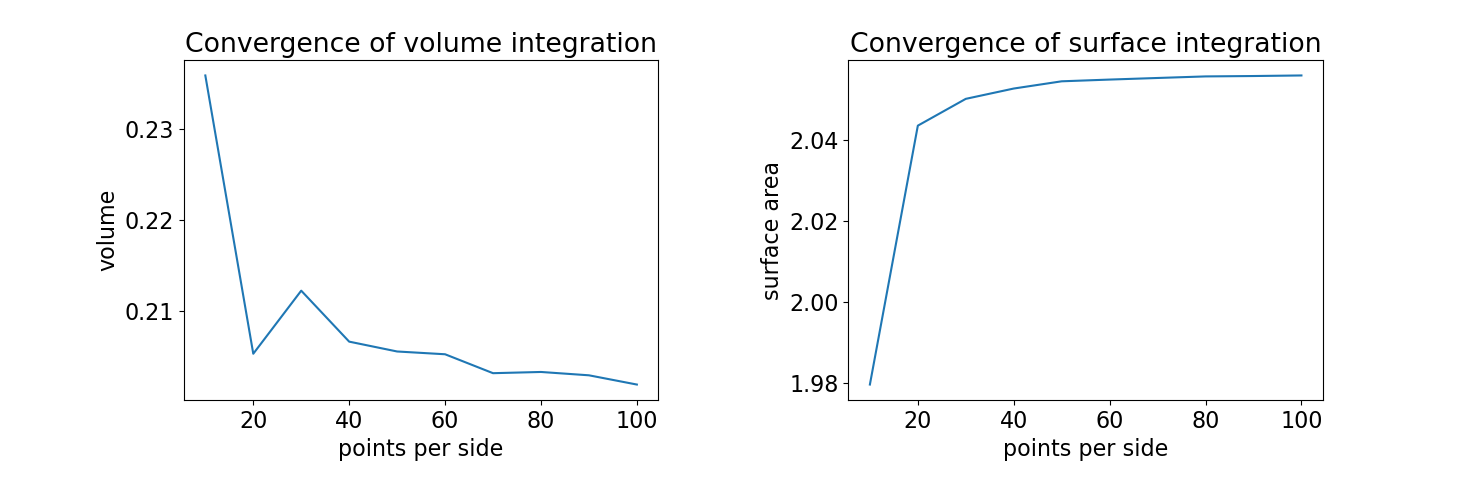}
\caption{We verify that the proposed integration grid is sufficiently refined to obtain converged values of the volume and surface area of the geometry defined by the level set.}
\label{integration_convergence}
\end{figure}

The source term is defined through the manufactured solution and Eq. \eqref{kpp}. The solution field is discretized with a two-hidden-layer network of width $30$, which corresponds to $1080$ trainable parameters. We take the convergence criteria as $\mathcal{B}_f=1\times{10}^{-2}$, $\mathcal{Z}_f=5\times10^{-3}$, and $\mathcal{R}_f=1 \times 10^{-2}$. The learning rate for ADAM optimization is $5 \times 10^{-3}$. Unless otherwise noted, all further neural network discretizations will use hyperbolic tangent activation functions. The inhomogeneous Dirichlet BCs are defined by evaluating the assumed solution on the surface integration grid from marching cubes. We take a $75 \times 75 \times 75$ integration grid inside the unit cube to build the representation of the geometry. This corresponds to $82537$ interior integration points and $33160$ surface integration points. Figure \ref{integration_convergence} shows the convergence of the integral for the volume and surface area of the geometry with the given background integration grid and marching cubes triangulation of the surface. We observe that the given integration grid is sufficiently refined to accurately reproduce these two geometric quantities.

\begin{figure}[hbt!]
\centering
\includegraphics[width=1.0\textwidth]{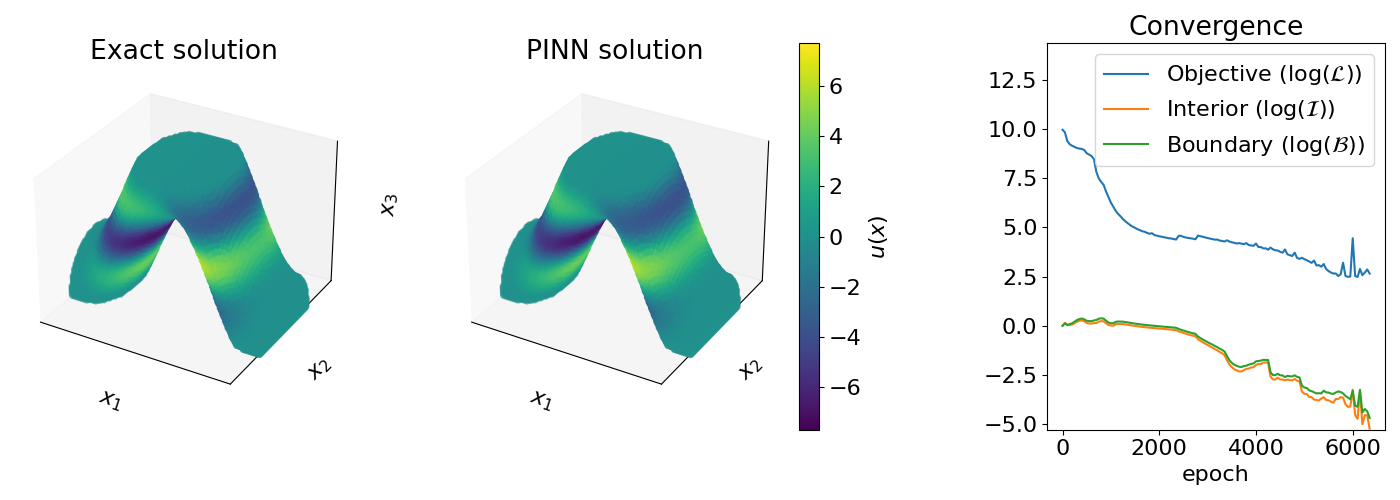}
\caption{Results of the PINN simulation of the Fisher-KPP equation defined on the two-way branching geometry with inhomogeneous Dirichlet BCs. The Augmented Lagrangian method is used to enforce the boundary conditions.}
\label{ex1}
\end{figure}

\begin{figure}[hbt!]
\centering
\includegraphics[width=1.0\textwidth]{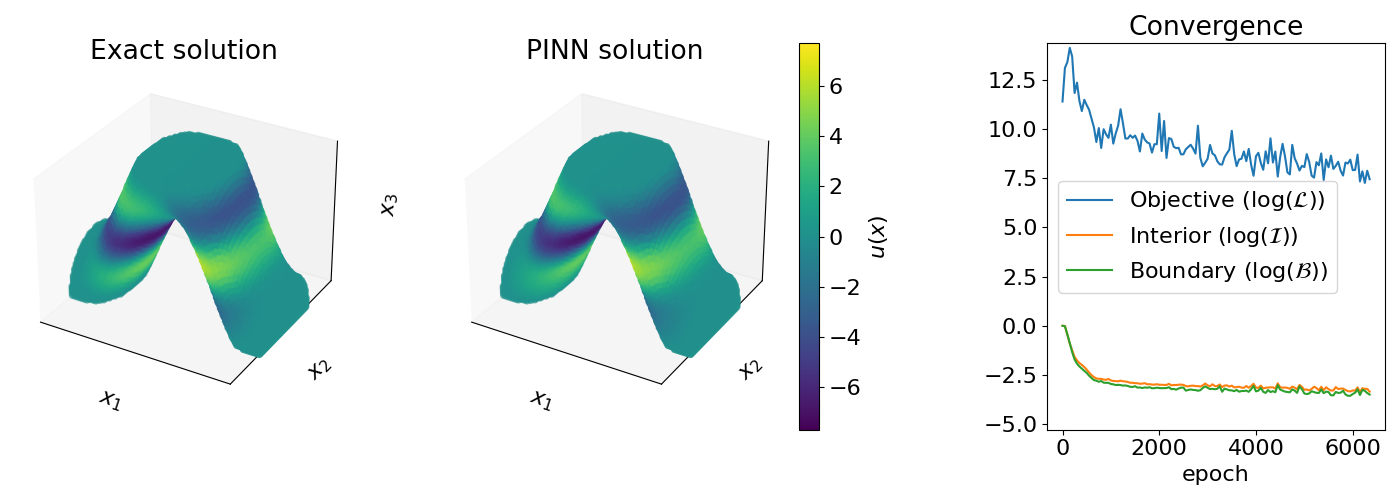}
\caption{Results of the PINN simulation of the Fisher-KPP equation on the branching geometry using the learning rate annealing method. We see that the interior and boundary error stagnate at large values compared to those obtained by the Augmented Lagrangian.}
\label{ex1_lra}
\end{figure}

\paragraph{} See Figure \ref{ex1} for the results. Using our proposed methodology, the interior and boundary errors at the converged solution are $\mathcal{B}=0.89\%$ and $\mathcal{I}=0.51\%$, respectively. Going forward, we will also report the relative $L_2$ error for the interior to facilitate comparisons with the literature. We call this error measure $\mathcal I_2$. The relative $L_2$ error with our proposed method on this problem is $\mathcal I_2=5.9 \times 10^{-3}$. The run time is $1269.8$ seconds for $6385$ optimization steps. Next, we run the learning rate annealing method with the same network and learning rate for the same number of steps to compare the performance. Per the recommendation of the authors \cite{wang_understanding_2020}, we set the moving average parameter in Eq. \eqref{lra} to $\alpha=0.9$. See Figure \ref{ex1_lra} for the results of this simulation. The interior and boundary errors after the allotted training period are $\mathcal{B}=3.6\%$ and $\mathcal{I}=2.9\%$, and the run time is $1297.8$ seconds. The relative $L_2$ error with the learning rate annealing method is $\mathcal I_2=3.9 \times 10^{-2}$. In this example, the Augmented Lagrangian achieves markedly better results.

\subsection{Heat transfer with nonlinear Robin boundary and tabletop geometry}

\paragraph{} The second test case will be another scalar-valued static problem with a nonlinearity in the BC rather than the differential operator. We consider steady-state heat transfer of a functionally graded three-dimensional geometry with a mix of Dirichlet and Robin BCs. The Robin BCs are a combination of a heat flux input and radiative heat loss. The problem is given by 
\begin{equation}
\begin{aligned}
    \nabla \cdot \Big( \kappa(\mathbf{x}) \nabla u(\mathbf{x}) \Big) + f(\mathbf{x}) = 0, \quad \mathbf{x} \in \omega,\\
    u(\mathbf{x}) = g(\mathbf{x}), \quad \mathbf{x} \in \partial \omega_D,\\
    -\nabla u (\mathbf{x}) \cdot \mathbf{\hat n} = q(\mathbf{x}) + \sigma u^4, \quad \mathbf{x} \in \partial \omega_R,
\end{aligned}
\end{equation}
where $\partial \omega_D$ is the Dirichlet region of the boundary, $\partial \omega_R$ is the Robin region, $\kappa(\mathbf x)$ is the thermal conductivity, $f(\mathbf x)$ is a volumetric heat source, and $q(\mathbf x)$ is a prescribed boundary heat flux. Though the Robin boundary is not explicitly discussed in Section \ref{sec:method}, it requires only modifying the Neumann section of the boundary to include the solution field in the BC along with the flux. The problem geometry is that of a four-legged ``tabletop'' defined by the following level set:
\begin{equation*}
\Phi(\mathbf{x}) = \Big( \cos 2\pi x - (1-5z^2)\Big)^2 + 9(y-0.5)^2 + \Big( \cos 2\pi y - (1-5z^2)\Big)^2 + 9(x-0.5)^2  - 3.
\end{equation*}

\begin{figure}[hbt!]
\centering
\includegraphics[trim = 0mm 50mm 0mm 35mm, clip, width=0.75\textwidth]{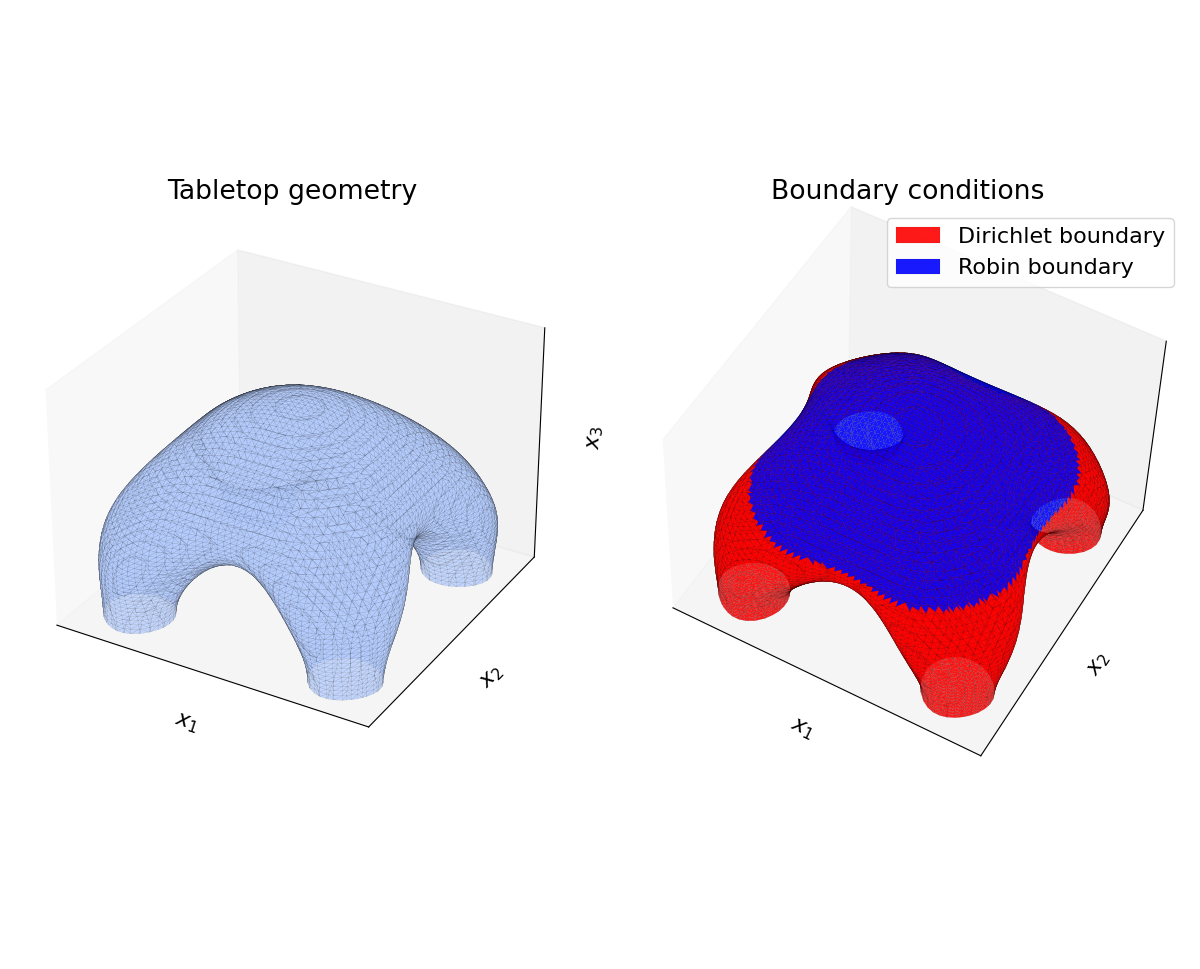}
\caption{Problem geometry for forced heat transfer with radiation BCs.}
\label{ex2_geo}
\end{figure}

The boundary is split into Dirichlet and Robin regions as follows:
\begin{equation*}
    \partial \omega_D = \{ \mathbf{x}| \Phi(\mathbf{x})=0, x_3 < 0.5\}, \quad \partial \omega_R = \{ \mathbf{x}| \Phi(\mathbf{x})=0, x_3 \geq 0.5\}.
\end{equation*}
See Figure \ref{ex2_geo} to visualize the geometry defined by the zero isocontour of this level set and the two regions of the boundary. Once again, the zero isocontour extends outside of the unit cube on which the background integration grid is defined. In this case, this means that the $x_3=0$ boundary surfaces are not detected by marching cubes. We enforce a homogeneous Dirichlet boundary along this surface in the discretization with 
\begin{equation*}
    \hat u(\mathbf{x} ; \boldsymbol \theta) = x_3\mathcal{N}(\mathbf{x} ; \boldsymbol \theta).
\end{equation*}

\begin{figure}[hbt!]
\centering
\includegraphics[width=0.975\textwidth]{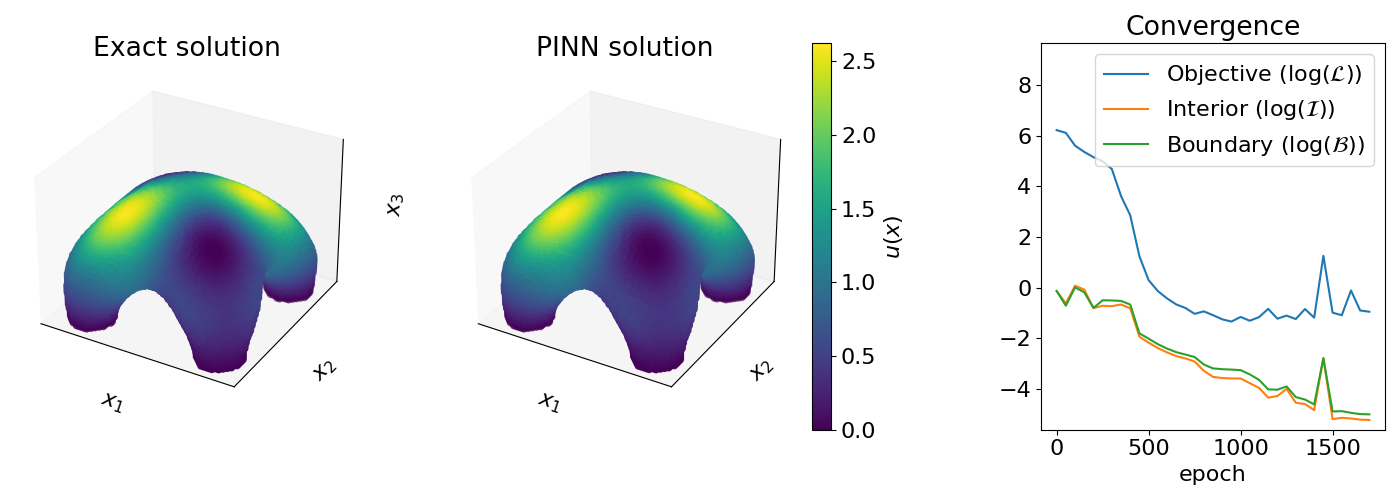}
\caption{Results of the heat conduction problem using the Augmented Lagrangian method on the tabletop geometry with a mix of inhomogeneous Dirichlet and radiative Robin BCs. The log of the absolute value of the objective function is shown in the convergence plot, given that it becomes negative due to the presence of the Lagrange multiplier term.}
\label{ex2}
\end{figure}

\begin{figure}[hbt!]
\centering
\includegraphics[width=0.975\textwidth]{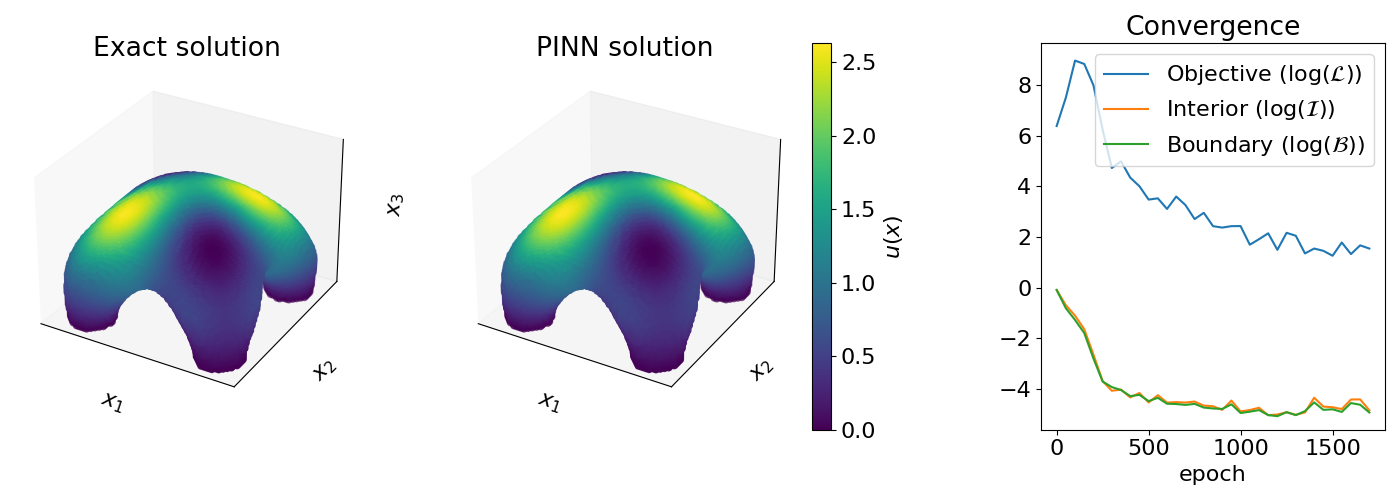}
\caption{Results of the heat conduction problem using the learning rate annealing method. In this example, there are no notable differences between the performance of the two techniques for BC enforcement.}
\label{ex2_lra}
\end{figure}

The forcing term and BCs will be computed through a manufactured solution. The assumed temperature field and the graded conductivity are given by 
\begin{equation*}
    u(\mathbf{x}) = 2x_3( 1 +  \sin(2\pi x_1)\sin(2\pi x_2)), \quad \kappa(\mathbf{x}) = 1 + x_3.
\end{equation*}
The forcing $f(\mathbf{x})$ is computed by substituting the manufactured solution into the governing equation for the bulk of the material. The applied heat flux $q(\mathbf{x})$ is computed at the boundary integration points with the manufactured solution per 
\begin{equation*}
    q(\mathbf{x}) = -\nabla u \cdot \mathbf{\hat n} - \sigma u^4, \quad \mathbf{x} \in \partial \omega_R.
\end{equation*}
We discretize the solution with a two-hidden-layer MLP with $30$ hidden units per layer, corresponding to $1080$ trainable parameters, and set the emissivity at $\sigma=0.1$. The learning rate for ADAM optimization is set at $5 \times 10^{-3}$. The convergence criteria are set at $\mathcal{B}_f=7.5\times{10}^{-3}$, $\mathcal{Z}_f=5\times10^{-3}$, and $\mathcal{R}_f=1 \times 10^{-2}$. The background integration grid is $75 \times 75 \times 75$, corresponding to $117486$ interior integration points and $47720$ surface integration points. In a manner equivalent to Figure \ref{integration_convergence}, we verify that this integration grid is sufficiently refined to obtain converged values of the volume and surface area. See Figure \ref{ex2} for results. The run time of this problem is $464.9$ seconds for $1702$ optimization steps. The two error measures at the converged solution are $\mathcal{B}=0.66\%$ and $\mathcal{I}=0.53\%$. The relative $L_2$ error is $\mathcal I_2=5.4 \times 10^{-3}$. The Augmented Lagrangian method converges in a small number of epochs even with a large number of nonlinear constraints on the boundary. We then run the learning rate annealing with the same network and learning rate, using an update weight of $\alpha=0.9$. We independently balance the gradients of the two loss terms corresponding to the Dirichlet and Robin boundaries per \cite{wang_understanding_2020}. See Figure \ref{ex2_lra} for the results of the simulation. After the allotted $1702$ steps, the learning rate annealing method obtains a boundary error of $\mathcal{B}=0.84\%$ and an interior error of $\mathcal I = 0.74 \%$ with a run time of $587.3$ seconds. The relative $L_2$ error on the interior is $\mathcal I_2=9.3 \times 10^{-3}$. We note that we loop through parameters to manually set gradients in with learning rate annealing as opposed to using a $\texttt{.backward()}$ call in PyTorch, which helps explain the increase in run time. In this example, the Augmented Lagrangian and learning rate annealing methods perform similarly.

\subsection{Linearly elastic non-prismatic pipe}

\paragraph{} We now perform a test on a vector-valued problem with inhomogeneous boundaries and a problem geometry with a different topology than the previous two examples. The problem is that of linear elastic deformation of a pipe with a varying cross-section. The governing equation for linear elasticity is 
\begin{equation}
\begin{aligned}
    \nabla \cdot \boldsymbol \sigma + \mathbf{f}(\mathbf{x}) = \mathbf{0}, \quad \mathbf{x} \in \omega, \\
    \mathbf{u}(\mathbf{x}) = \mathbf{g}(\mathbf{x}), \quad \mathbf{x} \in \partial \omega_D, \\
    \boldsymbol \sigma \cdot \mathbf{ \hat n} = \mathbf{t}( \mathbf{x}), \quad \mathbf{x} \in \partial \omega_N, \\
    \sigma_{ij} = C_{ijk\ell} \epsilon_{k\ell}, \quad \epsilon_{k\ell} = \frac{1}{2}\qty( \pd{u_k}{x_{\ell}} + \pd{u_{\ell}}{x_k}),
\end{aligned}
\end{equation}
\noindent where  $C_{ijk\ell}=\lambda \delta_{ij} \delta_{kl}+\mu(\delta_{ik} \delta_{jl} + \delta_{il} \delta_{jk})$ is the constitutive tensor corresponding to an isotropic, homogeneous linearly elastic material. The two Lame parameters are related to tabulated material parameters by
\begin{equation*}
    \lambda = \frac{E \nu}{(1+\nu)(1-2\nu)}, \quad \mu=\frac{E}{2(1+\nu)},
\end{equation*}
\noindent where $E$ is the Young's modulus and $\nu$ is the Poisson ratio. The non-prismatic pipe geometry is defined by the following level set:

\begin{equation*}
    \Phi(\mathbf{x}) = \Big( ( x_3 - 0.5)^2 + (x_2-0.5)^2 - R_2^2 - a^2\exp(-100(x_1-0.5)^2)  \Big) \Big( ( x_3 - 0.5)^2 + (x_2-0.5)^2 - R_1^2\Big),
\end{equation*}

\begin{figure}[hbt!]
\centering
\includegraphics[width=0.9\textwidth]{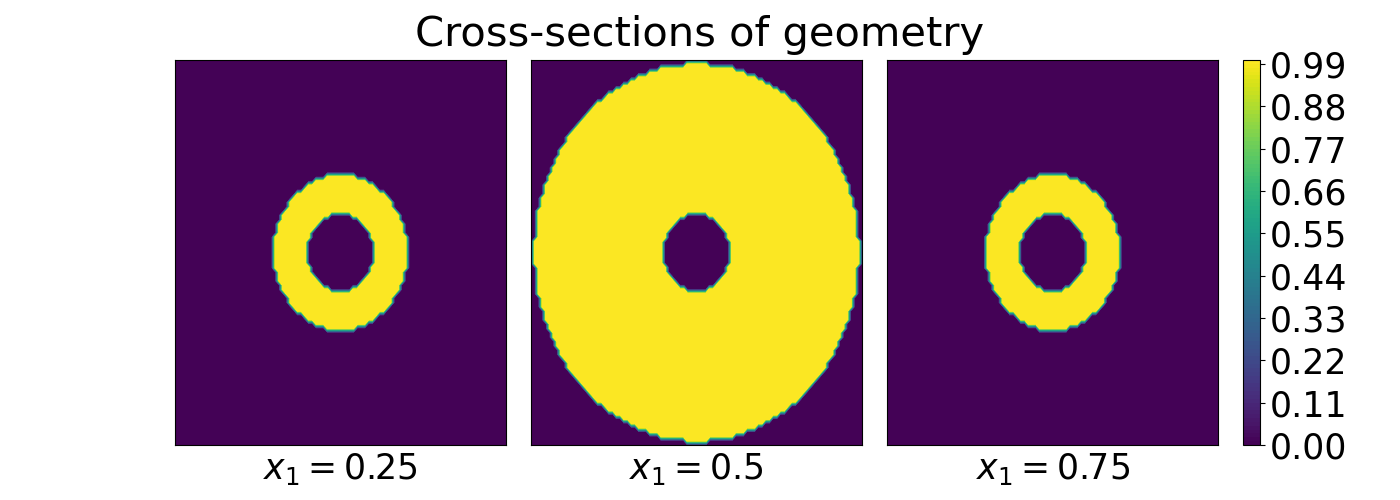}
\caption{The problem geometry given by the level set is that of a non-prismatic pipe, where the outer radius obtains a maximum at the center. In these plots, the inside of the volume is given by a phase value of $1$ and the outside by a phase value of $0$.}
\label{pipe}
\end{figure}

\noindent where $R_1$ and $R_2$ define the inner and outer radii of the pipe at the two edges of the boundary, and $a$ defines the increase of the outer radius at the center of the pipe. See Figure \ref{pipe} to visualize the geometry via cross-sections. We assume a manufactured solution corresponding to a radially symmetric expansion of the pipe whose magnitude is maximum at the center and falls off to zero at the two ends of the pipe. Defining the distance from the axis of the pipe as $r=\sqrt{(x_2-0.5)^2+(x_3-0.5)^2}$, the assumed solution is 
\begin{equation}\label{ex3man}
    \mathbf{u}(\mathbf{x}) = u_0\begin{bmatrix}
        0 \\ (x_2-0.5) r \sin(\pi x_1) \\
        (x_3-0.5) r \sin(\pi x_1)
    \end{bmatrix},
\end{equation}
\noindent where $u_0$ controls the magnitude of the displacement field. Note that this displacement field corresponds to clamped Dirichlet boundaries on the $x_1=0$ and $x_1=1$ faces of the domain. This is enforced automatically in the discretized solution with
\begin{equation*}
    \mathbf{\hat u}(\mathbf{x} ; \boldsymbol \theta) = \sin(\pi x_1)  \mathcal{N}(\mathbf{x};\boldsymbol \theta),
\end{equation*}
\noindent where $\mathcal{N}(\mathbf{x}; \boldsymbol \theta)$ is a vector-valued MLP neural network. The rest of the boundary will be inhomogeneous Dirichlet with values computed from evaluating the surface integration grid on the manufactured solution given by Eq. \eqref{ex3man}. Note that the error measures of Eqs. \eqref{error} are modified to accommodate the vector-valued solution field per 
\begin{equation*}
    \mathcal{I}(\boldsymbol \theta) = \frac{\int_{\omega} \lVert \mathbf{u}(\mathbf{x}) - \mathbf{\hat u}(\mathbf{x} ; \boldsymbol \theta)\rVert d\mathbf{x}}{\int_{\omega} \lVert \mathbf{u}(\mathbf{x})\rVert d\mathbf{x}}, \quad \mathcal{B}(\boldsymbol \theta) = \frac{\int_{\partial \omega} \lVert \mathbf{g}(\mathbf{x}) - \mathbf{\hat u}(\mathbf{x} ; \boldsymbol \theta) \rVert ds}{\int_{\partial \omega} \lVert \mathbf{g}(\mathbf{x}) \rVert ds}.
\end{equation*}

\begin{figure}[hbt!]
\centering
\includegraphics[width=0.99\textwidth]{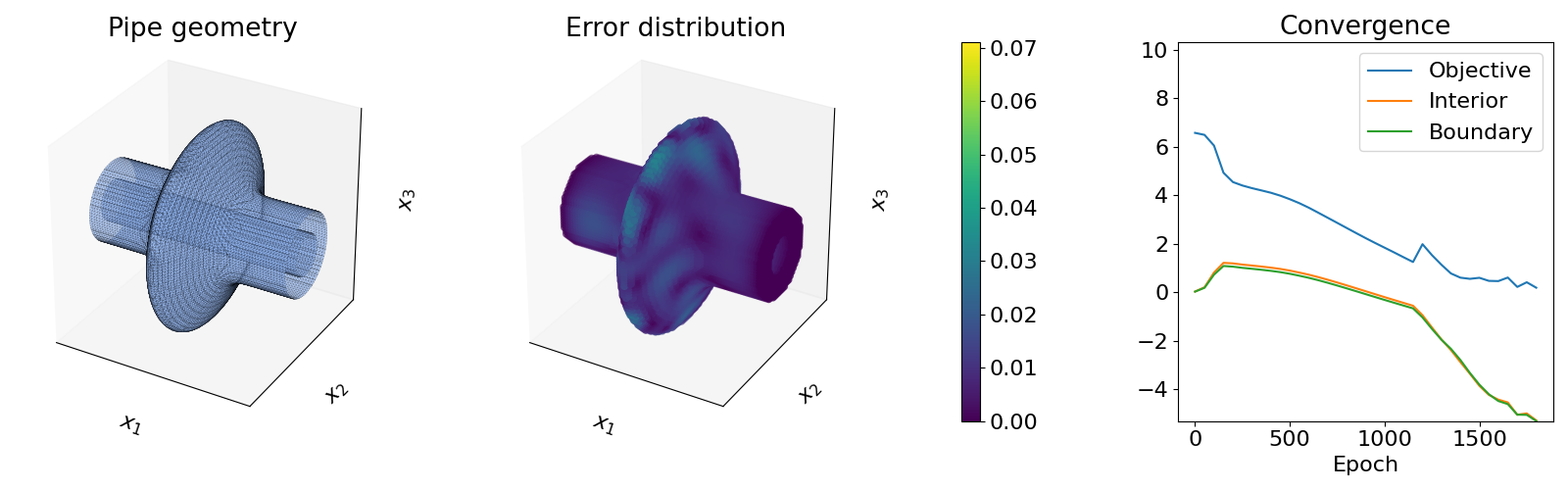}
\caption{Results of the PINN simulation of the non-prismatic pipe using the Augmented Lagrangian method to enforce inhomogeneous Dirichlet boundaries. Convergence is obtained in $1803$ optimization steps.}
\label{ex3results}
\end{figure}

\begin{figure}[hbt!]
\centering
\includegraphics[width=0.9\textwidth]{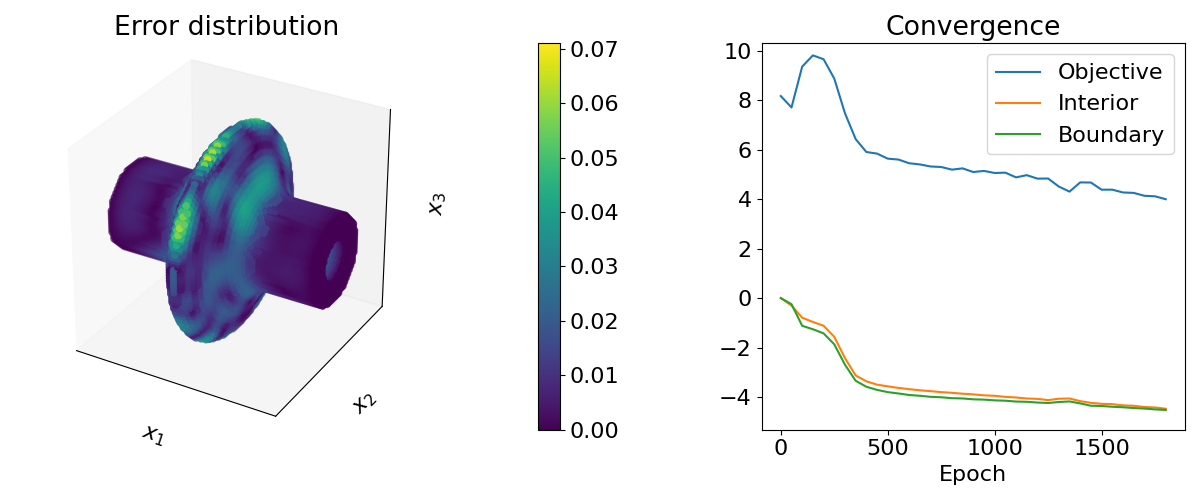}
\caption{Results of the PINN simulation of the non-prismatic pipe using the learning rate annealing method. We obtain larger errors when running the method for the same number of steps as the Augmented Lagrangian.}
\label{ex3_lra}
\end{figure}

The Augmented Lagrangian algorithm of Eq. \eqref{ALdis} is also modified to handle the vector-valued solution by replacing the squared residual of the strong form loss and the Dirichlet boundary with the squared magnitude of each residual. We take $u_0=25$ and use a two-hidden-layer neural network with $50$ hidden units per layer, corresponding to $2900$ trainable parameters, and a learning rate of $1 \times 10^{-3}$. The material properties are $E=1$ and $\nu=0.3$. The convergence criteria are set at $\mathcal{B}_f=5\times{10}^{-3}$, $\mathcal{Z}_f=2.5\times10^{-3}$, and $\mathcal{R}_f=1 \times 10^{-2}$, and the background integration grid is $75 \times 75 \times 75$ uniformly spaced points. This corresponds to $83164$ interior integration points and $45168$ surface integration points. We again verify that this integration grid is sufficiently refined to obtain converged values of the volume and surface area. See Figure \ref{ex3results} for a three-dimensional rendering of the geometry and the results of the simulation. The run time of this problem is $1855.8$ seconds. The two error measures at the converged solution are $\mathcal{I}=0.51 \%$ and $\mathcal{B}=0.49\%$. The relative $L_2$ error on the interior is $\mathcal I_2 = 4.2 \times 10^{-3}$. The method converges in $1803$ steps. We remark that in each of the three numerical examples, the inner loop problem in Algorithm \ref{alg:AL} converges no more than $11$ times, meaning that the maximum penalty parameter used in the Augmented Lagrangian method is $2^{11}=2048$. Finally, we re-run the simulation using the learning rate annealing method with the same network, learning rate, and number of steps as the Augmented Lagrangian. See Figure \ref{ex3_lra} for the results. The two error measures from learning rate annealing are $\mathcal{I}=\mathcal B = 1.1\%$, which are more than twice that of their counterparts from the Augmented Lagrangian. The relative $L_2$ error with the learning annealing method is $\mathcal I_2 = 1.2 \times 10^{-2}$. The run time of the numerical experiment is $1788.1$ seconds.

\section{Conclusion}
\label{sec:conclusion}

\paragraph{} The goal of this work was to discuss and study aspects of PINNs for PDEs to move them closer to a bona fide numerical solution framework, capable of applying to an equally large class of physical models, geometries, and BCs as the finite element method. To this end, we began with a discussion of the loss function, arguing that the strong form loss is the best choice if generality and exploiting the expressivity of the neural network are primary concerns. Next, we conducted an extensive review of techniques for BC enforcement. We ruled out coordinate transformations, distance functions, and PINN-FEM on first principles, as they struggle with Neumann boundary conditions and require significant pre-processing steps. We then investigated five different methods: the standard penalty approach, learning rate annealing, SA-PINNs, Lagrange multipliers, and the Augmented Lagrangian method. The results of a comparison of these methods with first-order optimization on a test problem are shown in Figure \ref{study}. The three best performing were SA-PINNs, the Augmented Lagrangian method, and learning rate annealing. Given that SA-PINNs dramatically increase the size of the optimization problem, we narrowed our focus down to the latter two methods. Finally, we found that Lagrange multipliers with quasi-Newton optimization methods struggled to handle the large number of nonlinear constraints arising from the neural network discretization. Thus, we ultimately settled on the Augmented Lagrangian and learning rate annealing as the two candidate techniques for BC enforcement, but chose to emphasize the Augmented Lagrangian in our exposition of a solution framework for three-dimensional geometries. We then studied the proposed method as well as learning rate annealing on three three-dimensional example problems with complex geometries and BCs. Using the method of manufactured solutions to obtain exact solutions, we showed that both techniques for BC enforcement can obtain accurate approximations with a modest parameter count and run times on the order of minutes. In the two-dimensional study, the Augmented Lagrangian outperformed learning rate annealing. This was also the case in two of the three three-dimensional numerical examples. While these results are not definitive, the Augmented Lagrangian method is the preferred strategy for handling BCs with PINNs in the opinion of the authors.

\paragraph{} In conducting a critical survey and head-to-head comparison of strategies for boundary condition enforcement, and in synthesizing existing strategies in the literature to enforce Dirichlet, Neumann, and/or Robin BCs on arbitrary three-dimensional domains, our work pushes PINNs one step closer to having the same impressive range of capabilities as incumbents such as FEM. However, many challenges remain. As seen in Figure \ref{conv}, we encountered convergence issues in the PINN problem when using a quasi-Newton optimizer with constraints. Though second-order optimization methods have recently become popular in the PINNs literature \cite{korbit_exact_2024, urban_unveiling_2025, kiyani_optimizing_2025}, there has been limited discussion of challenges unique to second-order methods in the context of non-convex problems. One example of this is the convergence of Newton methods to saddle points and maxima \cite{dauphin_identifying_2014}. The convergence issues we encountered and the tailoring of second-order methods to avoid concave regions of the loss landscape are opportunities for further research. Future work on this solution method will focus on extensions to dynamic PDEs on three-dimensional domains, as well as devising techniques to streamline the training of networks representing high-frequency solution fields. The numerical examples contained in this work show that standard PINN-based solutions are an accurate and efficient technique when the solution field is smooth and minimally oscillatory. However, numerical experimentation indicates that high-frequency solution fields present challenges for PINN solutions over and above the complexity of the geometry and/or BCs. See Appendix B for a preliminary discussion of using the Augmented Lagrangian for BC enforcement in the context of highly varying solution fields. The problems with oscillatory solutions have been noted, and various remedies have been proposed in the PINNs literature \cite{rahaman_spectral_2019, wang_eigenvector_2021, sitzmann_implicit_2020, ko_vs-pinn_2025}. Future work is required to explore the efficacy of these architectural modifications in the context of three-dimensional geometries and the Augmented Lagrangian method for BC enforcement. Devising methods that are accurate, reliable, and efficient on complex geometries and in the presence of high-frequency solutions will represent a decisive step in establishing PINNs as a mature numerical method. Finally, future studies might focus on the choice of optimizer with the Augmented Lagrangian method for BC enforcement. While ADAM performed satisfactorily on the problems we studied, we anticipate that more sophisticated optimization strategies may be required for solution fields exhibiting highly localized features and complex dynamics.

\section*{Acknowledgments}

\paragraph{} This work was funded by the National Defense Science and Engineering Graduate Fellowship (NDSEG) through the Department of Defense (DOD) and the Army Research Office (ARO). A. Doostan was supported by the Department of Energy, National Nuclear Security Administration, Predictive Science Academic Alliance Program (PSAAP) Award Number DE-NA0003962.

\section*{Supplementary material}

\paragraph{} The code used to generate the results and example problems presented in this work is available at https://github.com/crowan14/3D\textunderscore PINNs.

\appendix

\counterwithin*{equation}{section}
\renewcommand\theequation{\thesection\arabic{equation}}

\section{Distance functions for hard BC enforcement}

\paragraph{} In this appendix, we clarify and elaborate on the discussion in Section 3 regarding the use of distance functions for hard enforcement of Dirichlet and Neumann boundary conditions. In particular, we focus on issues related to the smoothness of distance functions and the increased computational cost of the method, putting aside concerns about increased complexity of the implementation and the necessity of offline computations. As discussed in \cite{sukumar_exact_2022}, distance functions obeying the Eikonal equation fail to be continuously differentiable when the boundary curve is not continuously differentiable. Intuitively, these distance functions have sharp ridges emanating from corners of the boundary curve. For example, an Eikonal distance function for the unit square $\mathbf x \in [0,1]^2$ is given by:

\begin{equation*}
    D(x_1,x_2) = \text{min}(x_1,x_2,1-x_1,1-x_2).
\end{equation*}

\begin{figure}[hbt!]
\centering
\includegraphics[width=0.65\textwidth]{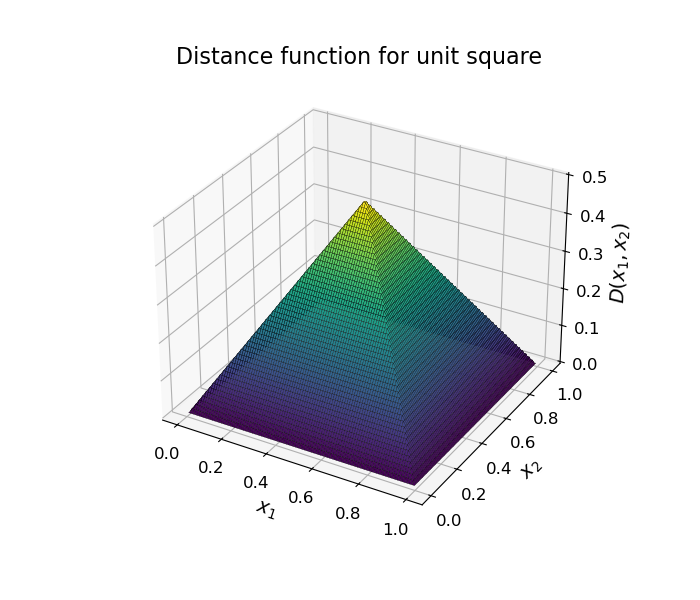}
\caption{The exact distance function for a unit square domain. The distance function is not continuously differentiable everywhere in the domain owing to the shape of the boundary curve.}
\label{dfunc}
\end{figure}

Recall that by an Eikonal distance function, we mean that $\lVert \nabla D\rVert=1$ and $D=0$ along the boundary. See Figure \ref{dfunc} for a visualization of the distance function for the unit square. It has a square pyramid shape, with four ridge line features with discontinuous derivatives. Because we rely on the strong form loss of the governing equation, we require two derivatives of the solution, and thus of the distance function, multiplying the neural network discretization. While the authors in \cite{sukumar_exact_2022} rely on approximate distance functions, even these have features that are not continuously differentiable near the boundary. Similar to an example given in \cite{goschel_enforcing_2025}, we illustrate the problems that arise from sharp features in the distance function with a numerical example. In fact, we can demonstrate the failure of training even in one spatial dimension, where the distance function has a derivative discontinuity only at a single point.

\paragraph{} The test problem we consider is a one-dimensional linear elliptic problem with homogeneous Dirichlet boundary conditions. The governing differential equation is 

\begin{equation*}
    \pdd{u}{x}+f(x)=0, \quad u(0)=u(1)=0,
\end{equation*}

\noindent where we enforce the two Dirichlet boundaries with a distance function approach. The exact distance function corresponding to the unit interval is $D(x)=0.5-|x-0.5|$. In this example, we use a smoothed absolute value to study the dependence of the numerical solution on the derivative discontinuity. In particular, our smoothed distance function is given by 

\begin{equation}\label{smooth_dist}
    D(x;\epsilon) = 0.5-(\sqrt{(x-0.5)^2+10^{\epsilon}} - 10^{\epsilon}).
\end{equation}

\begin{figure}[hbt!]
\centering
\includegraphics[width=0.95\textwidth]{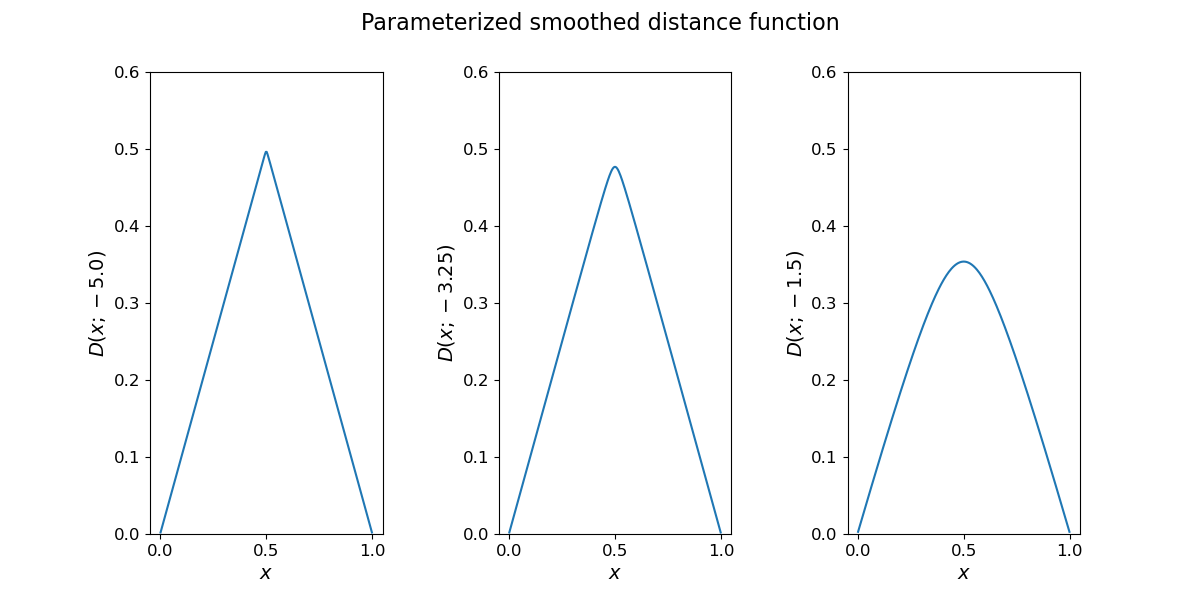}
\caption{The distance function we use to enforce homogeneous Dirichlet boundaries has adjustable smoothness properties. As $\epsilon \rightarrow -\infty$, we recover the true distance function, with a discontinuous derivative at the center of the domain.}
\label{smoothed}
\end{figure}

\begin{figure}[hbt!]
\centering
\includegraphics[width=0.99\textwidth]{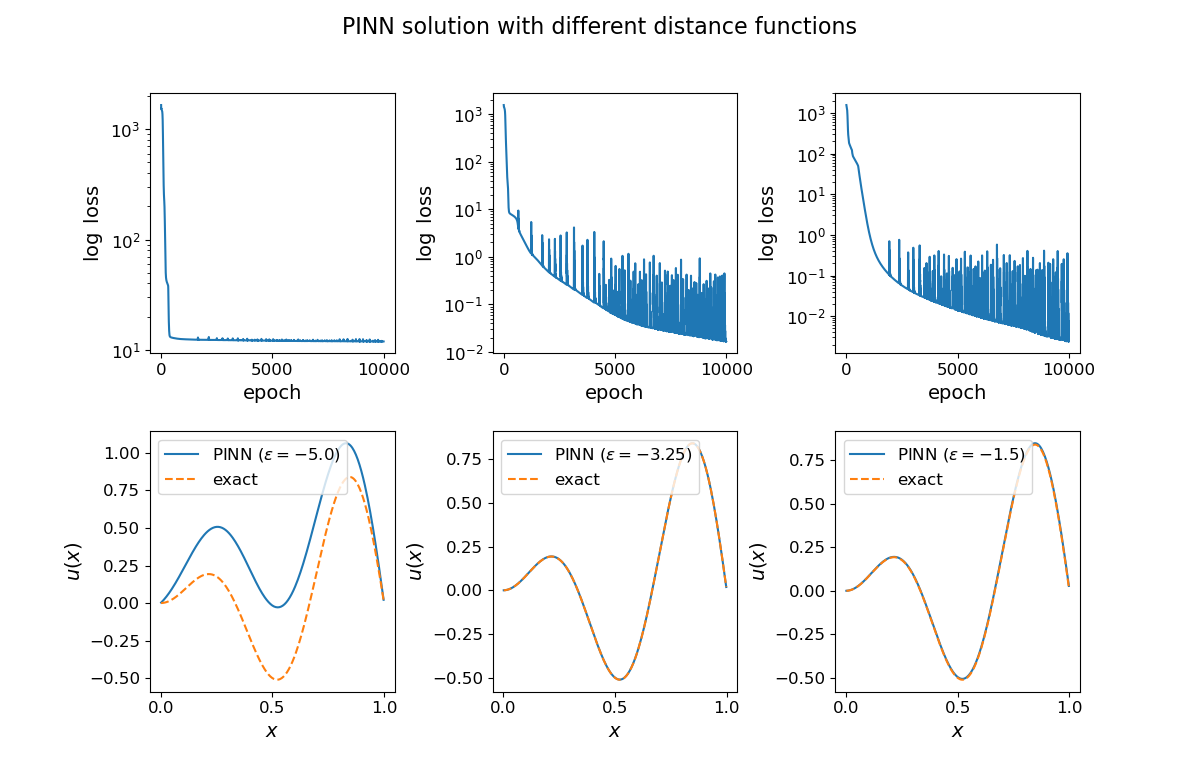}
\caption{When the distance function enforcing the Dirichlet boundaries has a nearly discontinuous derivative at a single point in the domain, the training stagnates at high loss, leading to a poor approximation of the exact solution.}
\label{comparison_results}
\end{figure}

When $\epsilon=-5$, this function is a close approximation of the exact distance function, with a near discontinuous derivative at $x=0.5$. At $\epsilon=-1.5$, we still have that $D(0)=D(1)=0$, but the derivative discontinuity is smoothed out. See Figure \ref{smoothed} for a visualization of this parameterized distance function. The solution is approximated with $u(x;\epsilon) = D(x;\epsilon)\mathcal N(x;\boldsymbol\theta)$, where $\mathcal N$ is a two-hidden-layer MLP neural network of width $25$. We use the method of manufactured solutions to obtain a ground truth solution to compare against, assuming that $u(x)=x\sin(3\pi x)$, with the source $f(x)$ computed accordingly. The training objective is the strong form loss, where there is no need to enforce boundary conditions with penalties or constraints. The integration grid consists of $250$ equally spaced points, the learning rate is $5 \times 10^{-3}$, and we run $10000$ epochs of ADAM optimization. We find a numerical solution for $\epsilon=-5,-3.25,-1.5$, which correspond to distance functions of different smoothness properties. See Figure \ref{comparison_results} for a comparison of the three different solutions. When $\epsilon=-5$, the derivative of the distance function is nearly discontinuous, and the optimizer stagnates at a high loss value, leading to a poor approximation of the true solution. The two smoother distance functions lead to satisfactory recovery of the exact solution. This example shows that even for a smoothed distance function, a rapid change in the derivative at a single point in the domain is sufficient to cause training to fail. Taking note of this, the authors in \cite{sukumar_exact_2022} avoid these derivative discontinuities in approximate distance functions by moving integration points away from the domain boundary, but we believe this strategy to be ad hoc, error-prone, and potentially infeasible when boundary effects are important to the PDE solution. 

\paragraph{} Next, we explore the increased computational cost of the approach outlined in \cite{sukumar_exact_2022} to strongly enforce Neumann boundary conditions. Again, we use a one-dimensional test problem for the sake of illustration. To isolate the study of computational cost from issues with differentiability, we work with the smoothed distance function of Eq. \eqref{smooth_dist} $D(x; \epsilon=-2.5)$. Note that training was successful for $\epsilon=-3.25$ and $\epsilon=-1.5$, so we expect no convergence issues with this choice of the smoothness parameter. Again, defining the computational domain as the unit interval, we wish to solve the following boundary value problem:

\begin{equation*}
    \pdd{u}{x} =0, \quad \pd{u}{x}(0)=\pd{u}{x}(1)=1.
\end{equation*}

Note that the exact solution to this problem can be written down by inspection as $u(x)=x+c$, where $c$ is an arbitrary constant. This issue of non-uniqueness will be addressed shortly. First, to strongly enforce the Neumann BCs, the solution can be discretized as follows: 

\begin{equation}\label{hard_neumann}
    u^{\text{DIST}}(x; \boldsymbol \theta_1 , \boldsymbol \theta_2) = \mathcal N_1(x; \boldsymbol \theta_1) - D \pd{D}{x} \pd{\mathcal N_1}{x}(x; \boldsymbol \theta_1) - D(-1+2x) + D^2 \mathcal N_2(x ; \boldsymbol \theta_2),
\end{equation}

\noindent where $D$ is the smoothed distance function, $\mathcal N_1$ and $\mathcal N_2$ are MLP neural networks, and the function $q(x)=-1+2x$ is chosen to enforce the Neumann boundary conditions. Note that the effectiveness of this discretization relies on $|\pd{D}{x}(0)|=|\pd{D}{x}(1)|=1$, a condition which we verify is accurately satisfied by our smoothed distance function. Note that the discretization now relies on two neural networks, both of which have two hidden-layers with width $25$. The optimization problem governing the numerical solution is

\begin{equation*}
    \underset{\boldsymbol \theta_1 , \boldsymbol \theta_2}{\text{argmin }} \frac{1}{2} \int_0^1 \qty( \pdd{u^{\text{DIST}}}{x})^2 dx + \lambda u^{\text{DIST}}(0)^2,
\end{equation*}

\noindent where $\lambda=100$ is a penalty parameter, and the penalty on the solution at the left end is used to remove the non-uniqueness of the elliptic problem with two Neumann boundaries. Note that the Neumann boundary conditions are satisfied by construction, but the form of the discretized solution is more complicated, including twice as many trainable parameters, and relying on derivatives of both the distance function and the neural network $\mathcal N_1$. This means that in the strong form loss, third spatial derivatives of the network are required. Similarly, differentiability issues with the distance function may be even further amplified by this formulation, given the presence of $\partial D / \partial x$ in the discretization of Eq. \eqref{hard_neumann}. For this problem, we work with a smoothed distance function and are thus not concerned with its differentiability. To assess the increased computational cost of derivative computations with automatic differentiation, we compare the above formulation with a penalty-based formulation of the Neumann BCs given by 

\begin{equation*}
    \underset{\boldsymbol \theta}{\text{argmin }} \frac{1}{2} \int_0^1 \qty(\pdd{u^{\text{PEN}}}{x})^2 dx + \lambda\qty( u^{\text{PEN}}(0)^2 + \pd{u^{\text{PEN}}}{x}(0)^2 + \pd{u^{\text{PEN}}}{x}(1)^2 ),
\end{equation*}

\noindent where the penalty-based solution is discretized with a two-hidden-layer network of width $25$, e.g., $u^{\text{PEN}}(x; \boldsymbol \theta) = \mathcal N(x ; \boldsymbol \theta)$. For both problems, we run $10000$ epochs of ADAM with a learning rate of $5\times 10^{-3}$ and the same integration grid of $250$ equally spaced points. We confirm that both obtain the exact solution $u(x)=x$ by the end of training. In both cases, we time the training process, with the only difference between the two methods being the formulation of the loss and the mathematical form of the solution. The penalty method runs for approximately $7$ seconds, whereas the exact enforcement method runs for approximately $27$ seconds. The same results were obtained each time for a total of $5$ runs. On this particular example problem, this means that the increased complexity of the solution discretization, involving derivatives of both the distance function and one of the networks, corresponds to a nearly four-fold increase in the computational cost of the method.

\section{Stress-testing the PINNs solver with oscillatory solution fields}

\paragraph{} In this Appendix, we stress test the Augmented Lagrangian method for BC enforcement in the presence of oscillatory solution fields with localized features. To easily visualize the manufactured solution, we restrict attention to a two-dimensional example problem, with a geometry defined by the following level set:

\begin{equation*}
    \Phi(\mathbf x) = \qty(x_1 - \frac{1}{2})^2 + \qty(x_2 - \frac{1}{2})^2 - \frac{1}{4},
\end{equation*}

\noindent which represents a circle with radius $R=1/2$ centered at $\mathbf x = (1/2,1/2)$. The governing equation is the Poisson problem with an inhomogeneous Dirichlet boundary condition:

\begin{equation*}
\begin{aligned}
    \nabla^2 u(\mathbf  x) + f(\mathbf x) = 0 , \quad \mathbf x \in \omega, \\
    u(\mathbf x) = g(\mathbf x), \quad \mathbf x \in \partial \omega.
\end{aligned}
\end{equation*}

In the three-dimensional examples discussed previously, the manufactured solutions have been relatively smooth and low-frequency, giving rise to smooth, low-frequency boundary conditions. In this example, we choose an oscillatory manufactured solution with large gradients inside the domain and especially near the boundary:

\begin{equation*}
    u(\mathbf x) = \tanh\qty(\frac{1}{4}\sin \Big( 6\pi (x_1-1/2)  \sin(2\pi(x_2-1/2)) \Big) ).
\end{equation*}

\begin{figure}[hbt!]
\centering
\includegraphics[width=0.99\textwidth]{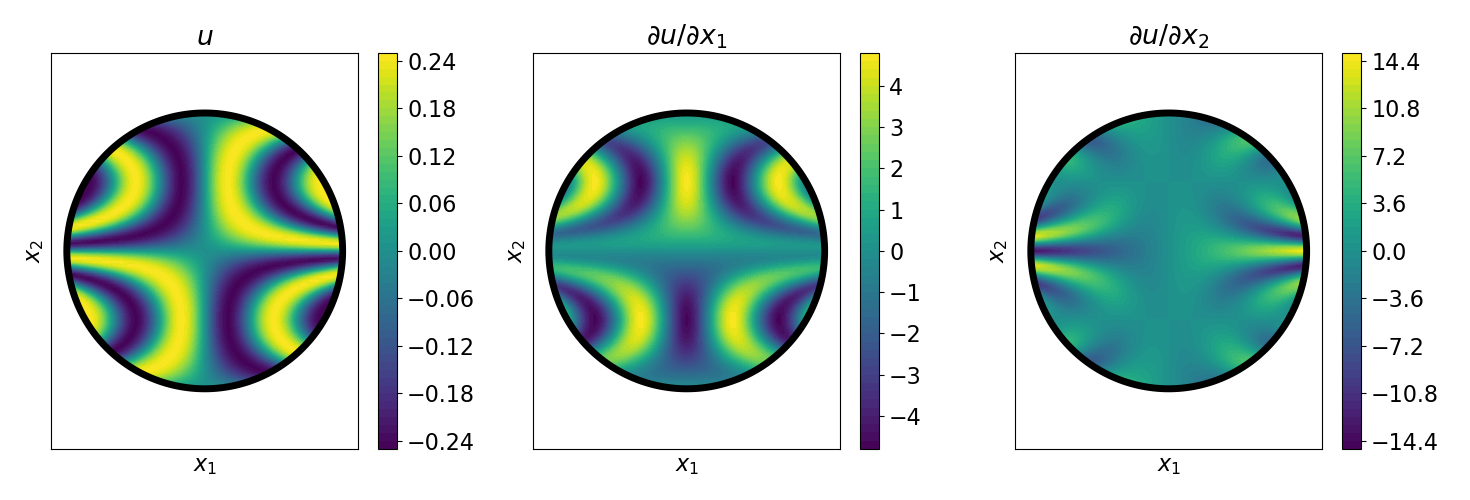}
\caption{The oscillatory manufactured solution and the two components of its gradient. There are especially small-scale features at the boundary along the left and right sides of the circular domain at $x_2=1/2$. The inhomogeneous Dirichlet boundary condition and the source term are computed using this assumed solution.}
\label{hf_man}
\end{figure}

\begin{figure}[hbt!]
\centering
\includegraphics[width=0.99\textwidth]{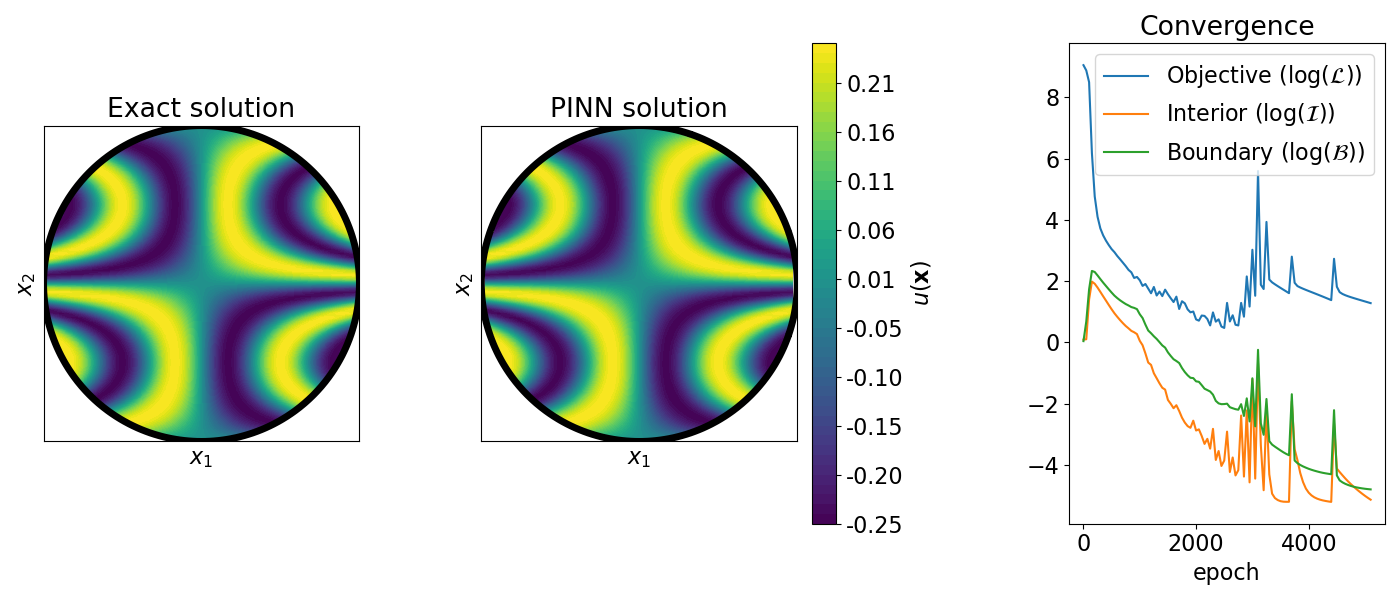}
\caption{PINN solution using Fourier feature embedding to resolve the oscillatory solution field. The Augmented Lagrangian method leads to accurate enforcement of the inhomogeneous Dirichlet boundaries when the network architecture is modified to account for the spectral bias.}
\label{hf_sol}
\end{figure}

See Figure \ref{hf_man} to visualize the manufactured solution and the two components of its gradient. The inhomogeneous Dirichlet boundary condition and source term $f(\mathbf x)$ are obtained from the manufactured solution. To compute integrals over the domain, we use a $75 \times 75$ uniform background grid on the unit square, which results in $4281$ points inside the circular domain. Surface integrals are computed using $200$ points along the boundary spaced at equal angular increments. Using ADAM optimization with a learning rate of $1 \times 10^{-3}$ and with a three-hidden-layer MLP network of width $50$, we do not obtain convergence with any settings of the convergence criteria $\mathcal B_f, \mathcal Z_f, \mathcal R_f$ we can find in an allotted $15000$ training epochs. We take this failure of convergence to be the consequence of the well-known challenges MLP networks have in representing oscillatory regression targets \cite{rahaman_spectral_2019}. Thus, in this case, the success of training hinges on the expressivity of the neural network discretization rather than the technique for BC enforcement. One popular remedy for the ``spectral bias'' of neural networks is the Fourier feature embedding method \cite{wang_eigenvector_2021}. With Fourier features, the spatial coordinate is mapped to a periodic feature $\boldsymbol \gamma$ with $\boldsymbol \gamma(\mathbf x) = [ \cos(2\pi \mathbf B \mathbf x) , \sin(2\pi \mathbf B \mathbf x)]^T$, where $B_{ij} \sim N(0,\sigma^2)$ is a random embedding matrix. The Fourier feature network is then an MLP applied to the Fourier feature vector $\boldsymbol \gamma(\mathbf x)$ rather than the raw spatial coordinate. Here, we provide an initial exploration of the Augmented Lagrangian for BC enforcement with oscillatory solution fields by using a Fourier network to solve the Poisson problem defined by the manufactured solution in Figure \ref{hf_man}. We take $10$ Fourier features with an embedding matrix whose entries are drawn from a standard normal distribution ($\sigma=1$). The subsequent MLP network has three hidden layers, each with a width $40$. The convergence criteria we use for the solution are $\mathcal R_f = 5 \times 10^{-4}$, $\mathcal Z_f=5 \times 10^{-3}$, and $\mathcal B_f = 1\times 10^{-2}$, where the gradient convergence criterion needs to be reduced with the Fourier feature network because the initial gradient is so large. See Figure \ref{hf_sol} for the results of the simulation. The two relative $L_1$ error metrics are $\mathcal I=0.56\%$ and $\mathcal B=0.82\%$, with a relative $L_2$ error in the interior of $\mathcal I_2=6.1 \times 10^{-3}$. This example demonstrates that with the appropriate choice of neural network architecture, the Augmented Lagrangian method leads to accurate solutions even on complex, oscillatory solution fields.

\end{document}